\documentclass[11pt]{article}

\usepackage{amssymb,latexsym,amsmath}

\usepackage{graphicx}

\begin{document}

\newcommand{\bfi}{\bfseries\itshape}

\makeatletter

\@addtoreset{figure}{section}

\def\thefigure{\thesection.\@arabic\c@figure}

\def\fps@figure{h, t}

\@addtoreset{table}{bsection}

\def\thetable{\thesection.\@arabic\c@table}

\def\fps@table{h, t}

\@addtoreset{equation}{section}

\def\theequation{\thesubsection.\arabic{equation}}

\makeatother

\newtheorem{thm}{Theorem}[section]

\newtheorem{prop}[thm]{Proposition}

\newtheorem{lema}[thm]{Lemma}

\newtheorem{cor}[thm]{Corollary}

\newtheorem{defi}[thm]{Definition}

\newtheorem{rk}[thm]{Remark}

\newtheorem{exempl}{Example}[section]

\newenvironment{exemplu}{\begin{exempl}  \em}{\hfill $\surd$

\end{exempl}}

\newcommand{\comment}[1]{\par\noindent{\raggedright\texttt{#1}

\par\marginpar{\textsc{Comment}}}}

\newcommand{\todo}[1]{\vspace{5 mm}\par \noindent \marginpar{\textsc{ToDo}}\framebox{\begin{minipage}[c]{0.95 \textwidth}

\tt #1 \end{minipage}}\vspace{5 mm}\par}

\newcommand{\ea}{\mbox{{\bf a}}}

\newcommand{\eu}{\mbox{{\bf u}}}

\newcommand{\ueu}{\underline{\eu}}

\newcommand{\ueo}{\overline{u}}

\newcommand{\oeu}{\overline{\eu}}

\newcommand{\ew}{\mbox{{\bf w}}}

\newcommand{\ef}{\mbox{{\bf f}}}

\newcommand{\eF}{\mbox{{\bf F}}}

\newcommand{\eC}{\mbox{{\bf C}}}

\newcommand{\en}{\mbox{{\bf n}}}

\newcommand{\eT}{\mbox{{\bf T}}}

\newcommand{\eL}{\mbox{{\bf L}}}

\newcommand{\eR}{\mbox{{\bf R}}}

\newcommand{\eV}{\mbox{{\bf V}}}

\newcommand{\eU}{\mbox{{\bf U}}}

\newcommand{\ev}{\mbox{{\bf v}}}

\newcommand{\eve}{\mbox{{\bf e}}}

\newcommand{\uev}{\underline{\ev}}

\newcommand{\eY}{\mbox{{\bf Y}}}

\newcommand{\eK}{\mbox{{\bf K}}}

\newcommand{\eP}{\mbox{{\bf P}}}

\newcommand{\eS}{\mbox{{\bf S}}}

\newcommand{\eJ}{\mbox{{\bf J}}}

\newcommand{\eB}{\mbox{{\bf B}}}

\newcommand{\eH}{\mbox{{\bf H}}}

\newcommand{\leb}{\mathcal{ L}^{n}}

\newcommand{\eI}{\mathcal{ I}}

\newcommand{\eE}{\mathcal{ E}}

\newcommand{\hen}{\mathcal{H}^{n-1}}

\newcommand{\eBV}{\mbox{{\bf BV}}}

\newcommand{\eA}{\mbox{{\bf A}}}

\newcommand{\eSBV}{\mbox{{\bf SBV}}}

\newcommand{\eBD}{\mbox{{\bf BD}}}

\newcommand{\eSBD}{\mbox{{\bf SBD}}}

\newcommand{\ecs}{\mbox{{\bf X}}}

\newcommand{\eg}{\mbox{{\bf g}}}

\newcommand{\paromega}{\partial \Omega}

\newcommand{\gau}{\Gamma_{u}}

\newcommand{\gaf}{\Gamma_{f}}

\newcommand{\sig}{{\bf \sigma}}

\newcommand{\gac}{\Gamma_{\mbox{{\bf c}}}}

\newcommand{\deu}{\dot{\eu}}

\newcommand{\dueu}{\underline{\deu}}

\newcommand{\dev}{\dot{\ev}}

\newcommand{\duev}{\underline{\dev}}

\newcommand{\weak}{\stackrel{w}{\approx}}

\newcommand{\mild}{\stackrel{m}{\approx}}

\newcommand{\strong}{\stackrel{s}{\approx}}

\newcommand{\weakdown}{\rightharpoondown}

\newcommand{\opg}{\stackrel{\mathfrak{g}}{\cdot}}

\newcommand{\opunu}{\stackrel{1}{\cdot}}
\newcommand{\opdoi}{\stackrel{2}{\cdot}}

\newcommand{\opn}{\stackrel{\mathfrak{n}}{\cdot}}
\newcommand{\opx}{\stackrel{x}{\cdot}}

\newcommand{\tr}{\ \mbox{tr}}

\newcommand{\Ad}{\ \mbox{Ad}}

\newcommand{\ad}{\ \mbox{ad}}

\renewcommand{\contentsname}{ }

\title{Sub-Riemannian geometry and Lie groups. Part II. Curvature of metric spaces, coadjoint orbits   and associated representations}

\author{Marius Buliga \\
 \\
IMB\\
B\^{a}timent MA \\
\'Ecole Polytechnique F\'ed\'erale de Lausanne\\
CH 1015 Lausanne, Switzerland\\
{\footnotesize Marius.Buliga@epfl.ch} \\
 \\
\and
and \\
 \\
Institute of Mathematics, Romanian Academy \\
P.O. BOX 1-764, RO 70700\\
Bucure\c sti, Romania\\
{\footnotesize Marius.Buliga@imar.ro}}

\date{This version: 05.07.2004}

\maketitle

\thispagestyle{empty}

\newpage

\thispagestyle{empty}

\tableofcontents

\section{Introduction}

\indent 

This paper is the third in a series \cite{buliga1}, \cite{buliga2}, dedicated to the fundamentals of 
sub-Riemannian geometry and its implications in Lie groups theory. We also bring to the attention of the reader the paper \cite{bulvodopis} for the more analytical aspects of the theory. 

The point of view of these papers is that what we are dealing with are manifestations  of the emerging  non-Euclidean analysis. This become visible especially in the metric study of Lie groups endowed with left invariant generating distributions. 

Metric profiles and their curvature first appeared in  \cite{buliga3}. They  give new insights into the behaviour of metric spaces. The subject is in close link with sub-Riemannian Lie groups because curvatures could be classified by comparison with (metric profiles of)  homogeneous spaces. The homogeneous spaces we are interested in  can be seen as  factor spaces of Lie groups with left invariant distributions.

The curvature of a metric space in a point  is, by definition, the rectifiability class of the metric profile 
associated to the point. There are two problems here that we are trying to solve. 

The first problem comes from the fact that to a metric space with some geometric structure we can associate several metric profiles. The most interesting are the profile as a metric space and the dilatation 
profile, associated with the dilatation structure of the space. 

The dilatation structure is the basic object in the study of differential properties of metric spaces of a certain type. (For example any Riemannian or sub-Riemannian manifold, with or without conical singularities, has a dilatation structure. But there are metric spaces  which are not admitting metric tangent spaces, but they admit dilatation structures).    This structure   tells us what is the 
good notion of analysis on that space. There is an infinite class of different such analysis and the classical one, which we call Euclidean, is only one of them.

To each metric profile of the space we  associate a curvature. We shall have therefore a metric curvature and a dilatation curvature. The metric curvature is more difficult to describe and in fact the dilatation curvature contains more complex informations. The first problem related to curvature is: 
what curvature we measure?

The second problem of curvature is related to the fact that the definition of curvature as a rectifiability 
class is too abstract to work with, unless we have a "good" set of representatives for the rectifiability 
classes. Suppose that we have made such a choice of the "good" class of representatives. Then we shall say that the curvature of a space in a point is the representant of the rectifiability class. The second problem of curvature is: what is a "good" class of representatives? 

This is obviously a subjective matter; anyway in this paper this subjectivity is revealed and, once revealed, the proposed choices become less subjective. 

The goal of this paper is to show that curvature in the sense of rectifiability class of the dilatation profile 
can be classified using coadjoint orbits representations. We establish here a  bridge between:
\begin{enumerate}
\item[-]  curvature notion in spaces which admit metric tangent spaces at any point,  and 
\item[-] self-adjoint  representations of algebras naturally associated with the structure of the tangent space. 
\end{enumerate}

\thispagestyle{empty}

In our opinion this is the second link between sub-Riemannian geometry and quantum mechanics. The 
first link can  be uncovered from 
 Buliga \cite{buliga2} section 5 "Case of the Heisenberg group", as explained in the section \ref{lastsec} of the present paper. 
 
 An appendix concerning uniform and conical groups is added to the paper. More information about the subject can be found in  Buliga \cite{buliga2}, sections 3 and 4. 

\newpage

\section{Distances and metric profiles}

  The references for the first  subsection  are Gromov \cite{gromov}, chapter 3, Gromov \cite{gromovgr},  and Burago \& al. \cite{burago} section 7.4.  There are several definitions of distances between metric spaces. The very fertile idea of introducing such distances belongs to Gromov.

  \subsection{Distances between metric spaces}

In order to introduce the Hausdorff distance between metric spaces, recall
the Hausdorff distance between subsets of a metric space.

\begin{defi}
For any set $A \subset X$ of a metric space and any $\varepsilon > 0$ set
the $\varepsilon$ neighbourhood of $A$ to be
$$A_{\varepsilon} \ = \ \bigcup_{x \in A} B(x,\varepsilon)$$
The Hausdorff distance between $A,B \subset X$ is defined as
$$d_{H}^{X}(A,B) \ = \ \inf \left\{ \varepsilon > 0 \mbox{ : } A \subset B_{\varepsilon} \ , \ B \subset A_{\varepsilon} \right\}$$
\end{defi}

By considering all isometric embeddings of two metric spaces $X$, $Y$ into
an arbitrary metric space $Z$ we obtain the Hausdorff distance between $X$, $Y$ (Gromov \cite{gromov} definition 3.4).

\begin{defi}
The Hausdorff distance $d_{H}(X,Y)$ between metric spaces $X$ $Y$ is the infimum of the numbers
$$d_{H}^{Z}(f(X),g(Y))$$
for all isometric embeddings $f: X \rightarrow Z$, $g: Y \rightarrow Z$ in a
metric space $Z$.
\end{defi}

If $X$, $Y$ are compact then $d_{H}(X,Y) < + \infty$.

The Hausdorff distance between isometric spaces equals $0$. The converse is also true  in the class of compact metric spaces (Gromov {\it op. cit.} proposition 3.6).

Likewise one can think about a notion of distance between pointed metric spaces. A pointed metric space is a triple $(X,x,d)$, with $x \in X$. Gromov \cite{gromovgr}  introduced the distance between pointed metric spaces $(X,x,d_{X})$ and $(Y,y,d_{Y})$ to be the infimum of all 
$\varepsilon>0$ such that there is a distance $d$ on the disjoint sum $X \cup Y$, which extends the distances on  $X$ and $Y$, and moreover 
\begin{enumerate}
\item[-] $d(x,y) < \varepsilon$, 
\item[-] the ball  $B(x,\frac{1}{\varepsilon})$ in $X$ is contained in the $\varepsilon$ neighbourhood of 
$Y$, 
\item[-] the ball  $B(y,\frac{1}{\varepsilon})$ in $Y$ is contained in the $\varepsilon$ neighbourhood of 
$X$.
\end{enumerate}

Denote by $[X,x,d_{X}]$ the isometry class of $(X,x,d_{x})$, that is the class of spaces $(Y,y,d_{Y})$ such that it exists an isometry  $f:X \rightarrow Y$ with the property $f(x)=y$. 

The Gromov distance between isometry classes of pointed metric spaces is almost  a distance, in the sense that whenever two  of the spaces $[X,x,d_{X}]$, $[Y,y,d_{Y}]$, $[Z,z,d_{Z}]$ have diameter at most equal to 2, then the triangle inequality for this distance is true. 
We  shall use this distance and the induced  convergence for isometry classes of 
the form $[X,x,d_{X}]$, with $diam \ X \ \leq 2$.

\subsection{Metric profiles}

 The notion of metric profile was introduced in  Buliga \cite{buliga3}.
 
We shall denote by $CMS$ the set of isometry classes of pointed compact metric spaces. The distance on this set is the Gromov distance between (isometry classes of) pointed metric spaces and the topology 
is induced by this distance.

To any locally compact metric  space we can associate a metric profile. 
  
\begin{defi}
 The metric profile associated to the 
 locally metric space $(M,d)$ is  the assignment (for small enough $\varepsilon > 0$) 
$$(\varepsilon > 0 , \ x \in M) \ \mapsto \  \mathbb{P}^{m}(\varepsilon, x) = [\bar{B}(x,1) , \frac{1}{\varepsilon} d] \in CMS$$
\label{dmprof}
\end{defi}

We can define a notion of metric profile regardless to any distance.

\begin{defi}
A metric profile is a curve $\mathbb{P}:[0,a] \rightarrow CMS$ such 
that:
\begin{enumerate}
\item[(a)] it is continuous at $0$,
\item[(b)]for any $b \in [0,a]$ and fixed $\varepsilon \in (0,1]$ we have 
$$d_{GH} (\mathbb{P}(\varepsilon b), \mathbb{P}^{m}_{d_{b}}(\varepsilon,x))  \  = \ O(b)$$
\end{enumerate}
We used here the notation $\mathbb{P}(b) = [\bar{B}(x,1),d_{b}]$ and $\mathbb{P}^{m}_{d_{b}}(\varepsilon,x) = [\bar{B}(x,1),\frac{1}{\varepsilon}d_{b}]$.

The metric profile is nice if 
$$d_{GH} (\mathbb{P}(\varepsilon b), \mathbb{P}^{m}_{d_{b}}(\varepsilon,x))  \  = \ O(b \varepsilon)$$
\label{dprofile}
\end{defi}

The metric profile $\varepsilon \mapsto \mathbb{P}^{m}(\varepsilon, x)$ of a metric space $(M,d)$ for  a fixed $x \in M$ is a metric profile in the sense of the definition \ref{dprofile} if and only if the space 
$(M,d)$ admits a tangent cone. Indeed, a tangent cone $[V,v,d_{v}]$ exists if and only if the following 
limit 
$$[V,v,d_{v}] \ = \ \lim_{\varepsilon \rightarrow 0} \mathbb{P}^{m}(\varepsilon, x)$$ 
exists. In this case the metric profile $\mathbb{P}^{m}(\cdot , x)$ can be prolonged to $\varepsilon = 0$. The prolongation is a metric profile in the sense of definition \ref{dprofile}. Indeed, we have still to check 
the property (b). But this is trivial, because for any $\varepsilon, b >0$,  sufficiently small, we have 
$$\mathbb{P}^{m}(\varepsilon b, x) \ = \  \mathbb{P}^{m}_{d_{b}}(\varepsilon,x)$$
where  $d_{b} \ = \ (1/b) d$ and $\mathbb{P}^{m}_{d_{b}}(\varepsilon,x) = [\bar{B}(x,1),\frac{1}{\varepsilon}d_{b}]$.

Note that in the definition of a nice metric profile  is not stated that $\displaystyle \mathbb{P}(0) = \mathbb{P}^{m}_{d_{b}}(0)$. 

The metric profile of a Riemannian homogeneous space is just  a curve in the space 
$CMS$, continuous at $0$. Likewise, if we look at a homogeneous regular  sub-Riemannian manifold, the metric profile is not depending on points in the manifold. 

\section{Sub-Riemannian manifolds} 

Classical references to this subject are Bella\"{\i}che \cite{bell} and Gromov \cite{gromo}. The interested reader is advised to look also to the references of these papers.  

In the literature on sub-Riemannian manifolds not everything written can be trusted. The source of errors 
lies in the first notions and constructions, mostly in the fact that obvious properties connected to 
a Riemannian manifold are not true (or yet unproven) for a sub-Riemannian manifold. Of special 
importance is the difference between a normal frame and the frame induced by a privileged chart. Also, 
there are many things in the sub-Riemannian realm which have no correspondent in the Riemannian 
case. 

To close this introductory comments, let us remark that  for the informed reader should be clear that in the sub-Riemannian realm there are problems even with the notion of manifold. Everybody agrees to define a sub-Riemannian manifold as in the definition \ref{defsr}. However, a better name for the object 
in this definition would be "model of sub-Riemannian manifold". We can illustrate this situation by the 
following comparison: the notion of sub-Riemannian manifold given in the definition \ref{defsr} is to 
the real notion what the Poincar\'e disk is to the hyperbolic plane. Only that apparently nobody found the real notion yet.

Let $M$ be a connected manifold. A distribution (or horizontal bundle) is a subbundle 
$D$ of $M$. To any point $x \in M$ there is associated the vectorspace $D_{x} \subset 
T_{x}M$. 

Given the distribution $D$, a point $x \in M$ and a sufficiently small  open neighbourhood $x \in U \subset M$, one can define on $U$ a filtration of bundles
as follows. Define first the class of horizontal vectorfields on $U$:
$$\mathcal{X}^{1}(U,D) \ = \ \left\{ X \in \Gamma^{\infty}(TU) \mbox{ : }
\forall y \in U \ , \ X(y) \in D_{y} \right\}$$
Next, define inductively for all positive integers $k$:
$$ \mathcal{X}^{k+1} (U,D) \ = \ \mathcal{X}^{k}(U,D) \cup [ \mathcal{X}^{1}(U,D),
\mathcal{X}^{k}(U,D)]$$
Here $[ \cdot , \cdot ]$ denotes vectorfields bracket. We obtain therefore a filtration $\displaystyle \mathcal{X}^{k}(U,D) \subset \mathcal{X}^{k+1} (U,D)$.
Evaluate now this filtration at $x$:
$$V^{k}(x,U,D) \ = \ \left\{ X(x) \mbox{ : } X \in \mathcal{X}^{k}(U,D)\right\}$$
There are
$m(x)$, positive integer, and small enough $U$ such that $\displaystyle
V^{k}(x,U,D) = V^{k}(x,D)$ for all $k \geq m$ and
$$D_{x} =  V^{1}(x,D) \subset V^{2}(x,D) \subset ... \subset V^{m(x)}(x,D)$$
We equally have
$$ \nu_{1}(x) = \dim V^{1}(x,D) < \nu_{2}(x) = \dim V^{2}(x,D) < ... < n = \dim M$$
Generally $m(x)$, $\nu_{k}(x)$
may vary from a point to another.

The number $m(x)$ is called the step of the distribution at 
$x$.

\begin{defi}
The distribution $D$ is regular if $m(x)$, $\nu_{k}(x)$ are constant on the manifold $M$. 

The distribution is completely non-integrable if for any $x \in M$ we have $\displaystyle 
V^{m(x)} = T_{x}M$. 
\label{dreg}
\end{defi}

\begin{defi}
A sub-Riemannian (SR) manifold is a triple $(M,D, g)$, where $M$ is a
connected manifold, $D$ is a completely non-integrable distribution on $M$, and $g$ is a metric (Euclidean inner-product) on the distribution (or horizontal bundle)  $D$.
\label{defsr}
\end{defi}

A horizontal curve $c:[a,b] \rightarrow M$ is a curve which is almost everywhere derivable and for 
almost any $t \in [a,b]$ we have $$\dot{c}(t) \in D_{c(t)}$$ 
The class of horizontal curves will be denoted by $Hor(M,D)$. 

The lenght of  a horizontal curve is 
$$l(c) \ = \ \int_{a}^{b} \sqrt{g(c(t)) (\dot{c}(t), \dot{c}(t))} \mbox{ d}t$$
The length depends on the metric $g$.

The Carnot-Carath\'eodory (CC) distance associated to the sub-Riemannian manifold is the 
distance induced by the length $l$ of horizontal curves:
$$d(x,y) \ = \ \inf \left\{ l(c) \mbox{ : } c \in Hor(M,D) \
, \ c(a) = x \ , \  c(b) = y \right\} $$

The Chow theorem ensures the existence of a horizontal path linking any two sufficiently 
closed points, therefore the CC distance it at least locally finite. 

We shall work further only with regular sub-Riemannian manifolds, if not otherwise stated. 

\subsection{Normal frames and privileged charts}
    
Bella\"{\i}che introduced the concept of privileged chart around a point $x \in M$.  

Let $(x_{1}, ... , x_{n}) \mapsto \phi(x_{1}, ... , x_{n}) \in M$ be a chart of $M$ around $x$ (i.e. $x$ has coordinates $(0,....,0)$).  Denote by $X_{1}, ... , X_{n}$ the frame of vectorfields  associated to 
the coordinate chart.  The chart is called adapted at $x$  (or the frame is called adapted) if the 
following happens: $X_{1}(x), ... , X_{\nu_{1}}(x)$ forms a basis of $V^{1}(x)$, 
$X_{\nu_{1}+1}(x), ... , X_{\nu_{2}}(x)$ form a basis of $V^{2}(x)$, and so on.  

Suppose that the frame $X_{1}, ... , X_{n}$ is adapted at (x). The degree of $X_{i}$ at $x$ is then $k$ if 
$X_{i} \in  V^{k} \setminus V^{k-1}$.

\begin{defi}
A chart  is privileged around the point $x \in M$ if it is adapted at $x$ and for any $i = 1, ... , n$ the 
function 
$$ t \mapsto d(x, \phi( ... , t, ...)) $$ 
(with $t$ on the position $i$) is exactly of order  $deg \  X_{i}$  at $t=0$. 
\end{defi}

Privileged charts always exist, as proved by Bella\"{\i}che  \cite{bell} Theorem 4.15. As an example 
consider a Lie group $G$ endowed with a left invariant distribution. The distribution is completely non-integrable if it is generated by the left translation of a vector subspace $D$ of the algebra 
$\mathfrak{g} = T_{e}G$ which bracket generates the whole algebra $\mathfrak{g}$. Then the exponential map is a privileged chart at the identity $e \in G$, but generically not privileged at 
$x \not = e$. 

Let $X$ be a vectorfield on $M$ and $x \in M$. The degree of $X$ at $x$ is the order of the function 
$$ t \in \mathbb{R} \rightarrow \mathbb{R} \ \  ,  t \mapsto d(x, \exp(tX)(x))$$
and it is denoted by $deg_{x} X$.  The vectorfield is called regular in an open set $U \subset M$ if 
$deg_{x} X$ is constant for all $x \in U$. 

Consider now a frame of vectorfields $X_{1}, ... , X_{n}$ defined in $U \subset M$, open set. Let 
$x \in U$ and define: 
$$V_{i}(x) \ = \ span \ \left\{ X_{k}(x) \mbox{ : } deg_{x} \ X_{k}  = i  \right\}$$

\begin{defi}
A frame $X_{1}, ... , X_{n}$, defined in an open set $U \subset M$,  is normal if all vectorfields 
$X_{k}$ are regular in $U$ and moreover  at any $x \in U$ and any $i = 1, ... , n$ one has: 
 $$V^{i}(x) \ = \ V_{1}(x) + ... + V_{i}(x)$$ 
 (direct sum).
 \end{defi}

A normal frame transforms the filtration into a direct sum. Each tangent  space decomposes 
as a direct sum of vectorspaces $V_{i}$. Moreover, each space $V^{i}$ decomposes in a direct 
sum of spaces $V_{k}$ with $k \leq i$. 

Normal frames exist. Indeed, start with a frame $X_{1}, ... , X_{r}$ such that for any $x$ $X_{1}(x), ... , 
X_{r}(x)$ form a basis for $D(x)$.  Associate now to any word $a_{1} .... a_{q}$ with letters in the alphabet $1, ... ,r$ the multi-bracket 
$$[X_{a_{1}}, [ ... , X_{a_{q}}] ... ]$$ 
One can add,  in the lexicographic order, elements to the set $\left\{ X_{1}, ... , X_{r} \right\}$, until a normal frame is obtained. 

To a normal frame $X_{1}, ... , X_{n}$ and the point $x \in U$ one can associate a privileged chart. Inded, such a chart is defined by: 
\begin{equation}
(a_{1}, ... , a_{n}) \in R^{n} \equiv T_{x}M  \mapsto \phi_{x}(\sum a_{i} X_{i}(x)) = \exp\left( \sum_{i=1}^{n} a_{i} X_{i} \right)(x)
\label{privchart}
\end{equation}
Remark how the privileged chart changes with the base-point $x$.

The intrinsic dilatations associated to a normal  frame, in a point $x$, are defined via a choice of a 
privileged chart  based at $x$. In such a  chart 
$\phi$,  for any $\varepsilon > 0$ (sufficiently small if necessary)  the dilatation is defined by 
$$\delta_{\varepsilon} (x_{i})  = (\varepsilon^{deg \ i} x_{i} )$$
With the  use of the privileged charts (\ref{privchart}), for any $\varepsilon > 0$ (sufficiently small if necessary)  the dilatations are 
$$\delta^{x}_{\varepsilon} \left(\exp\left( \sum_{i=1}^{n} a_{i} X_{i} \right)(x)\right) \  = \  
\exp\left( \sum_{i=1}^{n} a_{i} \varepsilon^{deg X_{i}}  X_{i} \right)(x)$$

In terms of vectorfields, one can use an intrinsic dilatation associated to the normal frame. This dilatation transforms 
$X_{i}$ into $$\Delta_{\varepsilon} X_{i} = \varepsilon^{deg \ X_{i}} X_{i}$$ 

One can define then, as in Gromov \cite{gromo}, section 1.4,  or Vodop'yanov \cite{vodopis3}, deformed vectorfields with respect to fixed 
$x$ by 
$$\hat{X}^{x}_{i}(\varepsilon)(y) \ = \  \Delta_{\varepsilon}\left(\left(\delta^{x}_{\varepsilon}\right)^{-1}*  X_{i}\right)(y)$$

When $\varepsilon \rightarrow 0$ the vectorfields $\hat{X}^{x}_{i}(\varepsilon)$ converge uniformly 
to a vectorfield $X^{N,x}_{i}$, on small enough compact neighbourhoods of $x$.

The nilpotentization of the distribution with respect to the chosen normal frame, in the point $x$,  is then 
the bracket 
\begin{equation}
[X,Y]_{N}^{x} = \lim_{\varepsilon \rightarrow 0}  [\hat{X}^{x}(\varepsilon) , 
\hat{Y}^{x}(\varepsilon)]
\label{nilbra}
\end{equation}
We have the equality
$$[X_{i},X_{j}]_{N}^{x}(x) \ = \ [X^{N,x}_{i}, X^{N,x}_{j}](x)$$

It is generically false that there are privileged coordinates around an open set in $M$. We can state this as a theorem. 

\begin{thm}
Let $(M,D,g)$ be a regular sub-riemannian manifold of topological dimension $n$. There are  $ \emptyset \not = U \subset M$ an open subset and $\phi: B \subset \mathbb{R}^{n} \rightarrow U$ such that  $\phi$ is a privileged chart for any $x \in U$ if and only if $(M,D,g)$ is a Riemannian manifold. 
\end{thm}

\paragraph{Proof.} If $(M,D,g)$ is Riemannian then it is known that such privileged charts exist. We have to prove the converse. Suppose there is a map $\phi: B \rightarrow  U$, $B$ open set in $\mathbb{R}^{n}$, $\phi$ surjective, such that for any $x \in U$ $\phi$ is privileged. Consider the frame $X_{1}, ... , X_{n}$ of vectorfields tangent to coordinate lines induced by $\phi$. Then this is a normal 
frame in $U$. Moreover 
$$[X_{i}, X_{j}] \ = \ 0$$ 
for any $i,j = 1, ... ,n$ therefore the nilpotentization bracket (\ref{nilbra}) in any point $x \in U$ is equal to 
$0$. According to Mitchell \cite{mit} theorem 1 (in this paper theorem \ref{mite1} section \ref{mitchsection}) 
the tangent cone in the metric sense to $x$ is the Euclidean $\mathbb{R}^{n}$. But this implies that 
the manifold (which is supposed regular) is Riemannian. 
\quad $\blacksquare$

\subsection{The nilpotentization at a point}

The nilpotentization   defines at point $x$ a nilpotent Lie algebra structure. This in turn gives a nilpotent group operation on the tangent space at $x$ (in the classical sense), by the use of the Baker-Campbell-Hausdorff formula. Denote this group $N(x)$. By construction the group $N(x)$ and its Lie algebra 
are identical as sets (otherwise said the exponential is the identity). What we have is a vectorspace 
endowed with a Lie bracket and a multiplication operation, linked by the Baker-Campbell-Hausdorff formula. 

 Consider the distribution $ND(x)$ on the group $N(x)$ which is obtained by left translation in $N(x)$ of 
the vectorspace  $D_{x} \subset N(x)$ (inclusion of a space in the Lie algebra of the group $N(x)$, seen 
as the tangent space to the neutral element of the group). We equally transport the metric $g(x)$ by nilpotent left translations in $N(x)$ and we denote the metric on the bundle $ND(x)$ by $g_{N}(x)$. We obtain a regular sub-Riemannian manifold $(N(x), ND(x), g_{N}(x))$. The associated Carnot-Carath\'eodory distance is denoted by $d_{N}^{x}$. 

As a metric space $(N(x), d_{N}^{x})$ is a cone based in the neutral element of $N(x)$. It is therefore a good candidate for being the metric tangent cone to the space $(M,d)$ at $x$. 

Endowed with the group operation and dilatations, $N(x)$ is a Carnot group.

\begin{defi}
A Carnot (or stratified nilpotent) group is a
connected simply connected group $N$  with  a distinguished vectorspace
$V_{1}$ such that the Lie algebra of the group has the
direct sum decomposition:
$$n \ = \ \sum_{i=1}^{m} V_{i} \ , \ \ V_{i+1} \ = \ [V_{1},V_{i}]$$
The number $m$ is the step of the group. The number
$$Q \ = \ \sum_{i=1}^{m} i \ dim V_{i}$$
is called the homogeneous dimension of the group.
\label{dccgroup}
\end{defi}

Any Carnot group admits a one-parameter family of dilatations. For any
$\varepsilon > 0$, the associated dilatation is:
$$ x \ = \ \sum_{i=1}^{m} x_{i} \ \mapsto \ \delta_{\varepsilon} x \
= \ \sum_{i=1}^{m} \varepsilon^{i} x_{i}$$
Any such dilatation is a group morphism and a Lie algebra morphism.

We collect some important facts to be known about Carnot groups:

\begin{enumerate}
\item[(a)] The metric  topology and uniformity of $N$ are the same as Euclidean
topology and uniformity respective.
\item[(b)] The ball $B(0,r)$ looks roughly like the box
$\left\{ x \ = \ \sum_{i=1}^{m} x_{i} \ \mbox{ : }
\| x_{i} \| \leq r^{i} \right\}$.
\item[(c)] the Hausdorff measure $\mathcal{H}^{Q}$ is group
invariant and the Hausdorff dimension of a ball is $Q$.
\item[(d)] there is a one-parameter group of dilatations, where a
dilatation is an isomorphism $\delta_{\varepsilon}$ of $N$ which
transforms the distance $d$ in $\varepsilon d$.
\end{enumerate}

The last item is especially important, because it leads to the introduction of an intrinsic notion 
of derivability in a Carnot group. 

In Euclidean spaces, given $f: R^{n} \rightarrow R^{m}$ and
a fixed point $x \in R^{n}$, one considers the difference function:
$$X \in B(0,1) \subset R^{n} \  \mapsto \ \frac{f(x+ tX) - f(x)}{t} \in R^{m}$$
The convergence of the difference function as $t \rightarrow 0$ in
the uniform convergence gives rise to the concept of
differentiability in it's classical sense. The same convergence,
but in measure, leads to approximate differentiability. 
This and
another topologies might be considered (see Vodop'yanov
\cite{vodopis}, \cite{vodopis1}).

In the frame of Carnot groups the difference function can be written using only dilatations and the group operation. Indeed, for any function between Carnot groups
$f: G \rightarrow P$,  for  any fixed point $x \in G$ and $\varepsilon >0$  the finite difference function is defined by the formula:
$$X \in B(1) \subset G \  \mapsto \ \delta_{\varepsilon}^{-1} \left(f(x)^{-1}f\left(
x \delta_{\varepsilon}X\right) \right) \in P$$
In the expression of the finite difference function enters $\delta_{\varepsilon}^{-1}$ and $\delta_{\varepsilon}$, which are dilatations in $P$, respectively $G$.

Pansu's differentiability is obtained from uniform convergence of the difference
function when $\varepsilon \rightarrow 0$.

The derivative of a function $f: G \rightarrow P$ is linear in the sense
explained further.  For simplicity we shall consider only the case $G=P$. In this way we don't have to use a heavy notation for the dilatations.

\begin{defi}
Let $N$ be a Carnot group. The function
$F:N \rightarrow N$ is linear if
\begin{enumerate}
\item[(a)] $F$ is a {\it group} morphism,
\item[(b)] for any $\varepsilon > 0$ $F \circ \delta_{\varepsilon} \
= \ \delta_{\varepsilon} \circ F$.
\end{enumerate}
We shall denote by $HL(N)$ the group of invertible linear maps  of
$N$, called the  linear group of $N$.
\label{dlin}
\end{defi}

The condition (b) means that $F$, seen as an algebra morphism,
preserves the grading of $N$.

The definition of Pansu differentiability follows:

\begin{defi}
Let $f: N \rightarrow N$ and $x \in N$. We say that $f$ is
(Pansu) differentiable in the point $x$ if there is a linear
function $Df(x): N \rightarrow N$ such that
$$\sup \left\{ d(F_{\varepsilon}(y), Df(x)y) \ \mbox{ : } \ y \in B(0,1)
\right\}$$
converges to $0$ when $\varepsilon \rightarrow 0$. The functions $F_{\varepsilon}$
are the finite difference functions, defined by
$$F_{t} (y) \ = \ \delta_{t}^{-1} \left( f(x)^{-1} f(x
\delta_{t}y)\right)$$
\end{defi}
    
For  differentiability notions adapted to general sub-Riemannian manifolds the reader can consult Margulis, 
Mostow  \cite{marmos1} \cite{marmos2}, Vodop'yanov , Greshnov \cite{vodopis2},  \cite{vodopis3} or Buliga \cite{buliga2}.

\section{Deformations of sub-Riemannian manifolds}

One can use privileged charts or normal frames to define  several deformations of a sub-Riemannian manifold around a point. These deformations can be described as curves  in the metric space $CMS$ of isometry classes of pointed compact metric spaces, 
 with the Gromov-Hausdorff  distance. For the isometry class of the pointed metric space $(X,x,d)$ we shall 
 use the notation $[X,x,d]$ or $[X,d]$ when he marked point is obvious. We shall work only with spaces 
 $(X,x,d)$ such that $X = \bar{B}(x,1)$. 

Consider   a privileged chart around $x \in M$. With this chart come the associated   dilatations $\delta^{x}_{\varepsilon}$.

The dilatation flow $\delta^{x}_{\varepsilon}$ induces the following deformation: 
 let $(D_{\varepsilon}, g_{\varepsilon})$ be the pair  
distribution - metric on the distribution obtained by transport with $(\delta^{x}_{\varepsilon})^{-1}$, namely: 
$$D\delta^{x}_{\varepsilon}(y) D_{\varepsilon}(y) =  D(\delta^{x}_{\varepsilon}y)$$
$$g(\delta^{x}_{\varepsilon}y)(D\delta^{x}_{\varepsilon}(y)u,D\delta^{x}_{\varepsilon}(y)v) = 
g_{\varepsilon}(y)(u,v)$$
for any $u,v \in T_{x}M$.  The deformation associated is 
\begin{equation}
[D, g,\varepsilon] = [ \bar{B}(p,1), (D, g,\varepsilon)]
\label{1defor}
\end{equation}
where the notation $(D, g,\varepsilon)$ is used  for the CC distance in the 
sub-Riemannian manifold $(M,D_{\varepsilon}, g_{\varepsilon})$.  

A slightly different deformation induced by the dilatation flow $\delta^{x}_{\varepsilon}$ is given by 
\begin{equation}
[\delta^{x}, \varepsilon] \ = \ [\overline{B}_{d}(x,\varepsilon), (\delta^{x},\varepsilon)]
\label{2defor}
\end{equation}
where the distance $(\delta^{x}, \varepsilon)$ is given by
$$(\delta^{x}, \varepsilon) (\delta^{x}_{\varepsilon}y,\delta^{x}_{\varepsilon}z) = d(y,z)$$ 
and the ball $\overline{B}_{d}(x,\varepsilon)$ is taken with respect to the original distance $d$. 

It is not trivial to remark that there is no reason for the equality $[D,g,\varepsilon] = [\delta^{x},\varepsilon]$. 

Another deformation is associated to the dilatations $\Delta_{\varepsilon}$ and pairs normal frame - 
Riemannian metric $g$.   This is simply
\begin{equation}
[\Delta,\varepsilon] = [\bar{B}(x,1), (\Delta,\varepsilon)]
\label{3defor}
\end{equation}
where $(\Delta,\varepsilon)$ is the Riemannian distance induced by the Riemannian metric $g_{\varepsilon}$ given by: 
$$g_{\varepsilon}(y) (\Delta_{\varepsilon}X(y),\Delta_{\varepsilon}Y(y)) = g(y)(X(y),Y(y))$$
for any pair of vectorfields $X, Y$. 

\section{Meaning of Mitchell theorem 1}
\label{mitchsection}

A key result in sub-Riemannian geometry is Mitchell \cite{mit} theorem 1:

\begin{thm}
For a regular sub-Riemannian space $(M,D,g)$, the tangent cone of
$(M,d)$ at $x \in M$ exists and it is isometric to $(N(x), d^{x}_{N})$.
\label{mite1}
\end{thm}

Recall that the limit in the Gromov-Hausdorff sense is defined up to isometry.
This means it this case that $N(D)$ is a model for the tangent space at $x$ to
$(M,d_{CC})$. In the Riemannian case $D = TM$ and $N(D) = R^{n}$, as a group with
addition. 

This theorem tells us nothing about the tangent bundle. 

One can identify in the literature several proofs of Mitchell theorem 1. Exactly what is proven in each available variant of proof? The answer is: each proof basically shows that various deformations, such as the ones introduced previously,  are metric profiles which can be prolonged   to $0$. Each of this metric profiles are close in the GH distance  to the original metric profile of the CC distance. More precisely: 
 
 \begin{lema}
 Let $\mathbb{P}'(t)$ be any of the previously introduced deformations $[\delta, \varepsilon]$, $[\Delta, \varepsilon]$, $[D,g,\varepsilon]$. Then 
 $$d_{GH}(\mathbb{P}^{m}_{\varepsilon}, \mathbb{P}'(\varepsilon)) = O(\varepsilon)$$
 \label{telema}
 \end{lema}
 
 The proof of this lemma reduces to a control problem. In the case of the profile $[\delta,\varepsilon]$, 
 this is Mitchell \cite{mit} lemma 1.2. 
      
  Mitchell \cite{mit} and Bella\"{\i}che \cite{bell} theorem 5.21, proposition 5.22,  proved the following: 
\begin{thm}
The deformation $\varepsilon \mapsto [D,g,\varepsilon]$ can be prolonged by continuity to $\varepsilon = 0$ and  the prolongation is a  metric profile. We have  
$$[D,g,0]  = [\bar{B}(0,1), d_{N}]$$
\end{thm}
The corresponding  result of Gromov \cite{gromo} section 1.4B and Vodop'yanov \cite{vodopis3}  is: 
\begin{thm}
The deformation $\varepsilon \mapsto [\Delta,\varepsilon]$ can be prolonged by continuity to $\varepsilon = 0$ and  the prolongation is a  metric profile. We have  
$$[\Delta,0]  = [\bar{B}(0,1), d_{N}]$$
\end{thm}

Finally, a similar theorem for the metric profile $[\delta^{x}, \varepsilon]$ is true.

Any of these theorems imply the Mitchell theorem 1, with the use of the approximation lemma \ref{telema}. But in fact these theorems are different statements in terms of metric profiles.

\section{Tangent bundles and dilatation structures}

There are several ways to associate a tangent bundle to a metric measure  space
(Cheeger \cite{cheeger} for what seems to be in fact a cotangent bundle) or to a regular sub-Riemannian manifold (Margulis \& Mostow
\cite{marmos1}, \cite{marmos2}, Vodop'yanov, Greshnov \cite{vodopis3}). These bundles differs. For example the  Cheeger tangent bundle of a sub-Riemannian manifold  can be identified with the distribution $D$ and Margulis-Mostow bundle is  the same as the usual tangent bundle, but with the fiber isomorphic with $N(D)$, instead of $R^{n}$.

For several reasons none of these constructions is completely satisfying in the case of sub-Riemannian manifolds. This is discussed at length in the paper Buliga, Vodop'yanov \cite{bulvodopis}. We adopt 
the point of view of this paper. The interested reader could also consult Buliga \cite{buliga2}, sections 3 to 5, 
for  more developed  constructions of the tangent bundle for sub-Riemannian Lie groups. 

A very important remark is that any notion of tangent bundle comes with an associated notion of differentiability. 

We begin by describing the Vodop'yanov tangent bundle. After this we explain what a dilatation structure  
is and how it induces another (but related) notion of tangent bundle.  

\begin{defi}
 The tangent bundle in the sense of Vodop'yanov is an assignment 
$x \in M \mapsto (\mathcal{O}(x), d_{N}^{x}, \delta_{\varepsilon}^{x})$ where $\mathcal{O}(x)$ is n open neighbourhood of 
$x$, $d_{N}^{x}$ is a quasi-distance on $\mathcal{O}(x)$ such that for any $y,z \in \mathcal{O}(x)$ we have 
\begin{equation}
d_{N}^{x}(\delta_{\varepsilon}^{x} y , \delta_{\varepsilon}^{x} z) = \varepsilon d_{N}^{x}(y,z) \  \  \ \forall 
\varepsilon \in [0,1]
\label{cond1}
\end{equation}
and 
\begin{equation} 
\frac{1}{\varepsilon} \mid d_{CC}(\delta_{\varepsilon}^{x} y , \delta_{\varepsilon}^{x} z) - 
d_{N}^{x}(\delta_{\varepsilon}^{x} y , \delta_{\varepsilon}^{x} z) \mid  \rightarrow 0 
\label{cond2}
\end{equation}
as $\varepsilon \rightarrow 0$, uniformly with respect to $y,z \in K \subset \mathcal{O}(x)$, compact. 

There is a nilpotent group operation on $\mathcal{O}(x)$, denoted by $\opx$. The neutral element of $\mathcal{O}(x)$ is 
$x$. 
\label{dtanvodopis}
\end{defi}

The differentiability notion introduced by Vodop'yanov is explained further. 

\begin{defi}
A function $f: M \rightarrow M'$ is differentiable at $x \in M$ in the sense of Vodop'yanov  if there is 
a group morphism $Df(x) : \mathcal{O}(x) \rightarrow O(f(x))$ which commutes with dilatations with the property:  
for any $\varepsilon >0$ there is $\mu>0$ such that for any $y \in \mathcal{O}(x)$, if $d_{CC}(x,y) \leq \mu$ then 
$$d_{N}^{f(x)} (f(y), Df(x)(y)) \leq \varepsilon d_{N}^{x}(x,y)$$
\label{difvodop}
\end{defi}

A related notion is a dilatation structure. We advertise  once again the point of view that this is the central notion of interest in the study of the differentiability properties of a sub-Riemannian manifold.

\begin{defi}
A dilatation structure associated to $(X,d)$ is an assignment $x \in X \mapsto \delta^{x}_{\varepsilon}: 
\mathcal{O}(x) \rightarrow \mathcal{O}(x)$, for all $\varepsilon \in (0,1]$, where $\mathcal{O}(x)$ is an open contractible neighbourhood of $x$ and all $\delta_{\varepsilon}^{x}$ are invertible, such that: 
\begin{enumerate}
\item[(a)] for any $x \in X$ the map 
$$ \mathbb{P}^{\delta}(x)(\varepsilon) \ = \ [\overline{B}_{d}(x,\varepsilon),(\delta,\varepsilon)]$$ 
is a nice metric profile, 
\item[(b)] if we denote by $\mathbb{P}^{m}(x)$ the metric profile of the space $(X,d)$ at $x$ then 
$$\lim_{\varepsilon \rightarrow 0} d_{GH}(\mathbb{P}^{m}(x)(\varepsilon), \mathbb{P}^{\delta}(x)(\varepsilon)) \ = \ 0$$
\item[(c)] for any $x \in X$ and $y \in \mathcal{O}(x)$ the map 
$$\left(\delta_{\varepsilon}^{\delta_{\varepsilon}^{x}y}\right)^{-1} \circ \delta_{\varepsilon}^{x}$$ 
converges uniformly to a map, as $\varepsilon$ tends to $0$. The convergence is uniform with respect 
to $x$ in a compact set. 
\item[(d)]  for small enough $\alpha > 0$ we have the uniform limit 
$$\delta_{\alpha}^{\delta^{x}_{\varepsilon} u} \rightarrow \delta_{\alpha}^{x}$$ 
which is uniform with respect to $x$ and $u \in \mathcal{O}(x)$, both in compact sets. 
\item[(e)] $\delta_{\varepsilon}^{x}$ contracts $\mathcal{O}(x)$ to $x$, uniformly with respect to $x$ in compact sets. 
\end{enumerate}
\label{dax}
\end{defi}
This definition generalizes the axioms for "uniform groups" introduced in Buliga \cite{buliga2}, section 3.1. A shortened version of this section can be found in the  Appendix. 

The point (a) from the definition \ref{dax} means that $\delta_{a}^{x} \circ \delta_{b}^{x}$ is 
approximately equal to $\delta^{x}_{b} \circ \delta^{x}_{a}$, because both are almost equal to  
$\delta^{x}_{ab}$ when we measure with the distance $d$. 

The point (b) from the definition is equivalent with the following assertion: for any $x_{0}$ and for any  $\varepsilon > 0$ there is $\lambda > 0$ and there is $\mu_{0}>0$ such that for any 
$x,y \in \overline{B}(x_{0},\lambda)$ and for any $\mu \in (0,\mu_{0}]$ we have
\begin{equation}
\frac{1}{\mu} \mid d(\delta_{\mu}^{x_{0}} x , \delta_{\mu}^{x_{0}} y) - \mu d(x,y) \mid \ \leq \ 
\varepsilon ( 1 + \max \left\{ d(x_{0},x), d(x_{0},y) \right\})
\label{whatbmeans}
\end{equation}

The meaning of the point (c)  will be explained a bit later, where we shall give a motivation for the definition \ref{dax}.

Remark that if $(X,d)$ admits a dilatation structure, it does not imply that it has a metric tangent cone 
at any of its points. 

We shall be concerned with $(X,d)$ which admits a dilatation structure. 

Such metric spaces exist. 

\begin{thm}
Any regular sub-Riemannian manifold admits a dilatation structure. 
\label{dilsr}
\end{thm}

\paragraph{Proof.} Indeed, such a dilatation structure is induced by any normal frame. This is a consequence of Mitchell theorem 1, more specifically it comes from lemma \ref{telema}, the  mentioned 
Mitchell type theorem for the metric profile $[\delta^{x}, \varepsilon] = \mathbb{P}^{\delta}(x)$ and proposition 6.3, section 6, Buliga, Vodop'yanov \cite{bulvodopis}.  Alternatively, the theorem is a direct consequence of the existence of a tangent 
bundle in the sense of Vodop'yanov, Greshnov \cite{vodopis3} and cited proposition 6.3, section 6 \cite{bulvodopis}.
\quad $\blacksquare$

In the case of regular sub-Riemannian manifolds we know that more it is true, namely that indeed the metric profiles can be prolonged to $\varepsilon = 0$. Therefore a metric tangent space exists in any point. Moreover, the dilatations $\delta_{\varepsilon}^{x}$ behave as homotethies for the nilpotent distance $d_{N}^{x}$.

Let us give now a motivation for the notion of dilatation structure. The most misterious seems to be 
the point (c) of the definition \ref{dax}. We shall call a map 
$$\left(\delta_{\varepsilon}^{\delta_{\varepsilon}^{x}y}\right)^{-1} \circ \delta_{\varepsilon}^{x}$$ 
an approximate infinitesimal left translation. Such maps appear from the following construction.

We are in a regular sub-Riemannian manifold, endowed with a normal frame and a tangent bundle 
in the sense of Vodop'yanov. 

Take $x \in M$ and $u,v \in \mathcal{O}(x)$ then we know that for any $\varepsilon > 0$ there is $\mu_{0} > 0$ 
such that  for any $\mu \in (0,\mu_{0}]$  we have 
$$\frac{1}{\mu} \mid d(\delta_{\mu}^{x}u, \delta_{\mu}^{x} v) - \mu d_{N}^{x}(u,v) \mid \ \leq \ \varepsilon$$
For small enough $\mu_{0}$ we can see the picture from the point of view of $\delta_{\mu}^{x}u$, that is 
in the same hypothesis it is also true that 
$$\frac{1}{\mu} \mid d(x, \delta_{\mu}^{x} v) -  d_{N}^{\delta_{\mu}^{x}u}(x,\delta_{\mu}^{x} v) \mid \rightarrow 0$$
as $\mu \rightarrow 0$. This is equivalent with 
$$\mid \frac{1}{\mu} d(x, \delta_{\mu}^{x} v) - d_{N}^{\delta_{\mu}^{x}u}(\left( \delta_{\mu}^{\delta_{\mu}^{x}u}\right)^{-1} (x), \left( \delta_{\mu}^{\delta_{\mu}^{x}u}\right)^{-1} \circ 
\delta_{\mu}^{x} (v) \mid \rightarrow 0$$
as $\mu \rightarrow 0$. We recognize here the appearance of the approximate left infinitesimal  translation. 

In this particular case these translations converge, as explained in the proof of theorem \ref{dilsr}. 
 We can actually use this trick of changing the 
point of view to prove that the limits, called left nilpotent translations,  are isometries with respect to the nilpotentized distance. 

The meaning of point (c) in definition \ref{dax} is that approximate infinitesimal translations converge. 
This allows to construct a tangent bundle.

\begin{defi}
The virtual tangent bundle associated to a dilatation structure is the assignment  
$$x \in X \ \mapsto \ VT^{\delta}_{x}X  \ = \ \left\{ \lim_{\varepsilon \rightarrow 0} 
\left(\delta_{\varepsilon}^{\delta_{\varepsilon}^{x}y}\right)^{-1} \circ \delta_{\varepsilon}^{x} \ \mbox{ : } 
y \in \mathcal{O}(x) \right\}$$
\label{dtandelta}
\end{defi}

Fix $x$ and denote 
$$u*v \ = \  \lim_{\varepsilon \rightarrow 0} 
\left(\delta_{\varepsilon}^{\delta_{\varepsilon}^{x}u}\right)^{-1} \circ \delta_{\varepsilon}^{x} (v)$$
Then $u*v$ is not an operation, but it leads to an operation, if we think at $u*v$ as the left translation of 
$v$ by "$u^{-1}$". Define for every $L \in VT^{\delta}_{x}X$ $$\delta_{\varepsilon} . L = \delta_{\varepsilon}^{x} \circ L \circ \delta_{\varepsilon}^{-1}$$

\begin{thm}
For any $\varepsilon>0$ sufficiently small and any $L \in VT^{\delta}_{x}X$ we have 
$\delta_{\varepsilon} . L \in VT^{\delta}_{x}$. Moreover $VT^{\delta}_{x}X$ forms a conical uniform group 
in the sense of definition 3.3 Buliga \cite{buliga2}. 
\label{tcon}
\end{thm}
The proof is a transcription of the proof of proposition 3.4 {\it op.cit.}, reproduced in  Appendix as proposition \ref{here3.4}. 

In the case of sub-Riemannian spaces more it is true: consider the dilatation structure induced by 
a normal frame. Then $VT^{\delta}_{x}M$ is isomorphic as a group with the metric tangent space 
at $x$ cf. Buliga \cite{buliga2} section 3.2.

\section{Curvature of a nice metric profile}

We shall define further the notion of curvature associated with a given metric profile. It is nothing but the rectifiability class of the metric profile with respect to the Gromov-Hausdorff distance. According to tastes 
and needs, one may use other distances than the mentioned Gromov-Hausdorff one, thus obtaining other notions of curvature. 

\begin{defi}
Two nice metric profiles $\mathbb{P}_{1}$, $\mathbb{P}_{2}$  are 
equivalent if 
$$d(\mathbb(P_{1}(a), \mathbb{P}_{2}(a)) = o(a)$$ 
The curvature class of a metric profile $\mathbb{P}$ is the equivalence class of $\mathbb{P}$. 
\end{defi}

This notion of curvature is too general from a geometer point of view. That is why we shall restrict 
to smaller classes of curvatures. Notice however that this is a choice of what should be a "good" curvature. For example we want  familiar spaces to have "good" curvature. A minimalist 
point of view could also be considered. "Good" curvatures should be associated to simple objects. 

Natural candidates which satisfy both requirements are homogeneous spaces.  Indeed, the metric profile of a homogeneous space does not depend on the base point and the space is simple in the sense that it is a factor space of two groups. 

Nevertheless, this choice is not based on a mathematical argument. At the moment where the paper is 
written  people  just like homogeneous spaces. Other classes of spaces could be interesting as well, 
for example factor spaces of Hopf algebras.

Let us explain what classification of curvatures means according to the preceding discussion. Later in the paper we shall refine the definition of  curvature like this:  if the metric profile (or other naturally associated profile) of  a metric space $(X,d)$ at a point $x \in X$ is  equivalent with the metric profile 
(or other...) of a particular geometric object $G$ (for example a homogeneous space), then we shall say 
that $(X,d)$ has curvature at $x$ and {\bf the} curvature of $(X,d)$ at $x$ {\bf is} $G$. 

Classification of curvatures means classification of (metric profiles or other profiles of) geometric objects from the class that we like. 

Here we shall  classify the curvatures  by comparison with 
deformations of  homogeneous ensembles, which are generalizations of homogeneous spaces. Until then we shall present some motivations for the notion of curvature proposed in this paper. 

This way of defining the curvature seems to be interesting even in the Riemannian case. It is partially connected to Cartan method of finding differential invariants, only that here the constructions are purely metric or intrinsicaly geometric. We (almost) don't use differential geometry methods.

\section{Homogeneous spaces and motivation for curvature}

A homogeneous metric space has the same metric profile associated to any of its points. That makes homogeneous spaces good candidates for classifying curvatures.

A homogeneous metric space is a metric space $(X,d)$ such that for any point $x \in X$ we have 
$$X = \left\{ f(x) \mbox{ :  } f \in Isom(X,d) \right\}$$
where $Isom(X,d)$ denotes the group of isometries of $(X,d)$.  Fix a point $p \in X$ and denote further 
$$Isom_{p}(X,d) \ = \ \left\{ f \in Isom(X,d) \mbox{ :  } f(x) = x \right\}$$
the stabiliser of the point $p$. 

 The coset class $Isom(X,d)/Isom_{p}(X,d)$ 
is (locally) homeomorphic with $(X,d)$ by the  map 
$$\pi: Isom(X,d)/Isom_{p}(X,d) \rightarrow X$$
which associates to  $p \in X$  the  coset $f Isom_{p}(X,d)$. 

The inner action of $Isom_{p}(X,d)$ on $Isom(X,d)$ gives an action of $Isom_{p}(X,d)$ on the coset space $Isom(X,d)/Isom_{p}(X,d)$. This inner action is compatible with the action of $Isom_{p}(X,d)$ on 
$X$ in the sense: for any $h \in Isom_{p}(X,d)$ and for all $f \in Isom(X,d)$ we have 
$$ \pi(h f h^{-1} Isom_{p}(X,D)) = h (\pi( f Isom_{p}(X,d)))$$

What is important is the geometry of a homogeneous space in the neighbourhood of one point. 
That is why we shall look to Lie algebras, endowed with a metric and a one parameter group of dilatations (equivalently, with a given gradation). 

Let $G$ be a Lie group and $G_{0}$ a subgroup. We shall imagine that there is a regular sub-Riemannian manifold $(X,d)$ such that  $G = Isom(X,d)$, $G_{0} = Isom_{p}(X,d)$. We want to see which  are the conditions upon $G$, $G_{0}$ in order for this to be possible. 

Any right-invariant 
vectorfield on $G$ descends on a vectorfield on left cosets  $G/G_{0}$. In particular, if we endow 
$G$ with a right-invariant distribution, then $G/G_{0}$ is endowed with a distribution induced by the 
descent of any right invariant "horizontal" frame.  $G/G_{0}$ is not usually  a regular sub-Riemannian manifold. Look for example to the case: $G = H(1)$, the 3 dimensional Heisenberg group. Take   $G_{0}$ as the one parameter group generated 
by an element of the distribution. Then $G/G_{0}$ is the Grushin plane, which is not a regular sub-Riemannian manifold. 

Consider on $G$ the right invariant distribution 
$$D" = Lie \  G_{0} + D'$$
If $G = Isom(X,d)$ and  $ G_{0} = Isom_{p}(X,d)$ then it should happen that  $Lie \ G_{0} \cap D' = 0$ and $D'$ descends on the distribution $D$ on $G/G_{0}$. 

Moreover  the action of $G_{0}$ on $G/G_{0}$  
correspond to   the action of $Isom_{p}(X,d)$ on $(X,d)$, which  is (the descent of) the inner action.  
We are not wrong if we suppose that isometries preserve the distribution at $0$, which translates into the 
condition: for any $h \in G_{0}$ 
$$Ad_{h} D' \subset D'$$
With the information from the preceding sections it is visible that this condition comes from the assumption that  isometries fixing the point $p$ are derivable. We don't have a proof for this assumption.  It is a natural condition that we impose. See further for another aspect of this condition.

We know one more thing about the homogeneous metric space $(X,d)$: its tangent cone.  Consider 
on $G$ with given distribution $D"$ the dilatations $\delta_{\varepsilon}$ and a privileged right-invariant 
basis around the neutral element. The knowledge of the tangent cone implies the following: 
\begin{enumerate}
\item[(a)]  we know some relations in the algebra $Lie \  G$ (described further), 
\item[(b)] we know that for any $h \in G_{0}$ $Ad_{h} \in HL(G,D")$, that is $Ad_{h}$ commutes with 
dilatations $\delta_{\varepsilon}$. 
\end{enumerate}

The previous discussion motivates the introduction of the following object. 

Let us consider triples $\sigma = ([ \cdot , \cdot ], \delta , g)$ where 
\begin{enumerate}
\item[(homs - a)] $[ \cdot , \cdot ]$ is a Lie bracket of a Lie algebra $\mathfrak{g}$, 
\item[(homs - b)] $\delta$ is a one parameter group of transformations $\varepsilon > 0 \mapsto \delta_{\varepsilon}$ which are simultaneously diagonalizable. The eigenspaces of any 
$\delta_{\varepsilon}$ form a direct sum decomposition of $\mathfrak{g}$ as a vectorspace,
$$\mathfrak{g} \ = \ V_{1} + ... + V_{m}$$ 
For any $k = 1, ... , m$ and any $x \in V_{k}, \varepsilon > 0$ we have $\delta_{\varepsilon} x = \varepsilon^{k} x$. 
\item[(homs - c)] $g$ is an inner product on $V_{1}$. Moreover, we shall suppose that  $V_{1}$ decomposes in an orthogonal sum 
$$V_{1} \ = \ D_{0} + D$$ 
such that $g =0$ on $D_{0}$ and $g$ is strictly positive defined 
on $D$.  We shall denote by $p$ the dimension of $D$. 
\item[(homs - d)] $D_{0}$ is a Lie subalgebra of $(\mathfrak{g}, \left[ \cdot , \cdot \right])$. Moreover, there is an 
isomorphism of Lie algebras between $D_{0}$ and a Lie subalgebra of $\mathfrak{so}(p)$, so we shall 
see the elements of $D_{0}$ as antisymmetric matrices acting on $D$ and the Lie bracket between 
$x_{0} \in D_{0}$ and $x \in D$ is 
$$[x_{0}, x] \ = \ x_{0} x$$
\item[(homs - e)] The direct sum decomposition  of $\mathfrak{g}$ is compatible with the bracket in the following sense: define first the deformed bracket 
$$[u, v]_{\varepsilon} \ = \ \delta_{\varepsilon}^{-1} [\delta_{\varepsilon} u , \delta_{\varepsilon}v]$$
for any $u,v \in \mathfrak{g}$. 
When $\varepsilon \rightarrow 0$ the deformed bracket converges: the limit exists 
$$\lim_{\varepsilon} [ u, v]_{\varepsilon} \ = \ [u,v]_{N}$$
and $\mathfrak{g}$ is a direct sum between the abelian Lie algebra $(D_{0}, [ \cdot , \cdot]_{N})$ and 
the Carnot algebra $D + V_{2} + ... + V_{m}$ with the bracket $[ \cdot , \cdot ]_{N}$ and dilatations 
$\delta_{\varepsilon}$.
\item[(homs - f)] finally, a differentiability condition: for any $x_{0} \in D_{0}$, $\varepsilon > 0$  and $y \in \mathfrak{g}$ we have: 
$$[x_{0}, \delta_{\varepsilon} x] \ + \ D_{0} \ = \ \delta_{\varepsilon} [x_{0}, x] \ + \ D_{0}$$
\end{enumerate}

\begin{defi}
We call homogeneous space any triple $\sigma$ which satisfies all  conditions (homs - a), ... , (homs - f).
\label{defhoms}
\end{defi}

\begin{rk}
This notion of homogeneous space is different from the usual one. We have not been able to find a special name for these triples. Notice however that in the realm of regular sub-Riemannian manifolds the previous definition introduces the natural notion of homogeneous space. 
\end{rk}

A morphism between triples $\sigma = (\mathfrak{g}, \delta, g)$ and $\sigma'=(\mathfrak{g}', \delta', g')$ 
is a linear map $F: \mathfrak{g} \rightarrow \mathfrak{g}'$ which is a Lie algebra morphism, it transforms  $g$ into $g'$ and it commutes with dilatations : $F \delta = \delta' F$.

Remark that instead of doing a factorisation we choose to introduce $D_{0}$ with null metric on it. 
This is equivalent to a metric factorisation.

 To a triple $\sigma = (\mathfrak{g}, \delta, g)$ is associated the length measure given by the metric $g$, translated at left on a neighbourhood of $0 \in \mathfrak{g}$ using the Baker-Campbell-Hausdorff formula. This 
length measure induces a (pseudo) distance denoted by $d$. For any $\varepsilon > 0$ define the (pseudo) distance 
$$d_{\varepsilon}(x,y) \ = \ d (\delta_{\varepsilon}^{-1} x , \delta_{\varepsilon}^{-1} y )$$
Make the metric factorisation of the domain of convergence of the Baker-Campbell-Hausdorff formula with respect to $d$. The resulting metric space is (in a neighbourhood of $0$)  a regular sub-Riemannian manifold. 

Declare $a>0$ admissible if the open ball centered at $0$ and radius $1$, defined with respect to the 
distance $d_{2a}$, is contained in (the factor of) the domain of convergence of the Baker-Campbell-Hausdorff formula.  
Consider now the number  
$$R = \min \left( 1, \sup \left\{a>0  \ \mbox{ admissible } \right\} \right)$$  
Associate then to $\sigma$ the compact metric space $(X_{\sigma}, d_{\sigma})$ where $X_{\sigma}$ is 
the closed ball of radius $1$ centered at $0$ and $d_{\sigma} = d_{R}$.  The function evoked previously  
is $$\sigma \mapsto [X_{\sigma}, d_{\sigma}, 0]$$

To $\sigma$ is associated also a metric profile, given by 
$$\mathbb{P}_{\sigma}(a)  = [\overline{B}(0, 1), d_{\sigma}(a), 0]$$ 
which is defined as  the deformation $[D, g,\varepsilon]$ described at (\ref{1defor}).  We call this the dilatation metric profile of $\sigma$.

From the point of view of this paper the rectifiability class of such a deformation is a good curvature. 

We give further a motivation for the fact that we consider such curvatures good.

\subsection{Riemannian surfaces}
 \label{rie}
 
 We shall begin by looking at homogeneous spaces which correspond to  Riemannian homogeneous $n$ manifolds. We are interested in triples $\sigma$ such that: 
 \begin{enumerate} 
 \item[-] the Lie algebra $\mathfrak{g}$ admits a decomposition 
$$\mathfrak{g}  \  = \ D_{0} + D$$
with dimension of $D$ equal to $n$, 
\item[-] On $D$ we have an Euclidean metric $g$, 
\item[-]  $D_{0}$ is a Lie subalgebra of $(\mathfrak{g}, \left[ \cdot , \cdot \right])$. Moreover, there is an 
isomorphism of Lie algebras between $D_{0}$ and a Lie subalgebra of $\mathfrak{so}(n)$, so we shall 
see the elements of $D_{0}$ as antisymmetric matrices acting on $D$ and the Lie bracket between 
$x_{0} \in D_{0}$ and $x \in D$ is 
$$[x_{0}, x] \ = \ x_{0} x$$
This limits the dimension of $D_{0}$ to be at most $n (n-1)/2$. 
\item[-] For any $\varepsilon > 0$ we have the (usual) dilatation $\delta_{\varepsilon} x \ = \ \varepsilon x$. 
\end{enumerate}

When $n=2$ we are looking to Riemannian 2 dimensional homogeneous surfaces. In this case  
$\mathfrak{g}$ is 2 or  3 dimensional. 

We start with the 3 dimensional case. Consider  a basis 
$X_{0}, X_{1}, X_{2}$ for $\mathfrak{g}$, such that $X_{0}$ generates $D_{0}$ and $X_{1}, X_{2}$ 
forming an orthonormal basis of $D$. The bracket relations that we know are: 
$$[X_{0}, X_{1}] =  a X_{2}$$
$$[X_{0}, X_{2}] = -a X_{1} $$
$$[X_{1}, X_{2}] = b X_{0} + c X_{1} + d X_{2} $$ 
We suppose that $a \not = 0$. 
From Jacobi identity we get $c=d= 0$. Therefore we have
$$[X_{0}, X_{1}] =  a X_{2}$$
$$[X_{0}, X_{2}] = -a X_{1}$$
$$[X_{1}, X_{2}] = b X_{0}$$
 We have a one dimensional family of homogeneous Riemannian surfaces, where the curvature can 
 be measured by $\mid ab \mid$. All these are 2 dimensional spheres of radius $1/\mid ab \mid$.  
 
 The second case is $dim \ D_{0} = 0$ and $\mathfrak{g}$ is a  2 dimensional Lie algebra. Take as previously a basis for $D = \mathfrak{g}$. We have only one bracket relation: 
 $$[X_{1}, X_{2}] = a X_{1} + b X_{2}$$
  This corresponds to a 2 dimensional family of homogeneous surfaces with negative curvature. When we look to deformations $[D,g ,  \varepsilon]$ and their equivalence as  metric profiles, we see that we are left with one dimensional family of negative curved surfaces, with curvature $\displaystyle -\sqrt{a^{2} +   b^{2}}$. 
  
  Therefore in the case of 2 dimensional riemannian manifolds, the curvature in the sense of this paper 
  is the Gauss curvature.

\subsection{About Cartan method of equivalence for sub-Riemannian manifolds}

For a description of  the Cartan method of equivalence for sub-Riemannian spaces consult Montgomery 
\cite{montgome}, chapter 7, especially sections 7.2, 7.8 (Riemannian surfaces), 7.7 (Distributions: torsion equals curvature), 7.10(Subriemannian contact three-manifolds). In the last mentioned 
section the work  of Hughen \cite{hughen} is followed in the exposition. 

We are not going to describe here the Cartan method of equivalence. The mentioned references should satisfy the interested reader. Instead, we shall comment the method. 

The main object of the method is the $G$-structure which encodes the sub-Riemannian manifold. Let us 
describe it. 

A coframe $(\eta_{1} , ... , \eta_{n})$ in $T^{*}M$ encodes  a regular sub-Riemannian manifold $(M,D,g)$ if we have 
$$D \ = \ ker \ (\eta_{1}^{p+1}, .... , \eta^{n}) \ \ \ , \ \ \ g \ = \ \left( (\eta_{1})^{2} + ... + (\eta_{n})^{2}\right) \mid_{D}$$
Such a coframe is called adapted to the sub-Riemannian manifold $(M,D,g)$. 

On can see the family of all adapted coframes to a regular sub-Riemannian manifold as a bundle over the manifold $M$,  with the structure structure group 
$$G_{0} \ = \ \left\{ \left( \begin{array}{cc} 
A & B \\
0 & C 
\end{array} \right)  \mbox{ :  } A \in O(p), \ B \in \mathcal{M}(p, n-p) ,  \ \ C \in GL(n-p, \mathbb{R}) \right\}$$
Such a bundle $\mathcal{B} \rightarrow M$ is called a $G_{0}$ structure. 

Cartan method of equivalence is a machine for solving the problem of local equivalence of $G_{0}$-structures. From this point of view, consider  two regular sub-Riemannian manifolds $(M,D, g)$ and 
$(M',D',g')$ with distinguished points $x \in M$ and $x' \in M'$. If there are open neighbourhoods 
$x \in U \subset M$ and $x' \in U' \subset M'$ and a diffeomorphism $\phi : U \rightarrow U'$ which induces a bundle isomorphism of the $G_{0}$ structures associated to $U$ and $U'$ then we shall 
say that $M$ and $M'$ are locally equivalent around $(x,x')$. 

We might be interested in  metric local equivalence between the pointed metric spaces $(M,x,d)$ 
and $(M',x',d')$, which means: there are open neighbourhoods 
$x \in U \subset M$ and $x' \in U' \subset M'$ and an isometry map $\phi : U \rightarrow U'$.

In this paper we are working with an even weaker notion of equivalence, that of metric profiles associated to deformations of sub-Riemannian manifolds. 

The author does not know 
any rigorous result concerning the relations between these equivalences. The reader under the temptation to use Rademacher or Stepanov type results for sub-Riemannian manifolds, like 
the one in Margulis, Mostow \cite{marmos2}, should be very careful.

The method described here as "classification" of curvatures is, at the present stage of development, 
only remotely connected to the modern applications of the Cartan method of equivalence. When applied 
to sub-Riemannian contact 3 manifolds, it seems to be closer to the work described in the paper 
 Agrachev, El Alaoui, Gauthier, Kupka \cite{aagk}, where the authors look for local normal form expansions for the exponential map associated to the metric sub-Riemannian structure.

\subsection{Homogeneous contact 3 manifolds}
\label{srcon}

  Contact manifolds are particular cases of sub-Riemannian manifolds. The contact distribution 
  is completely non-integrable. By using natural normalization of the contact form (see for example 
  Bieliavski, Flbel, Gorodski \cite{biefago} or Hughen \cite{hughen}) we can uniquely associate to a contact structure, endowed with a metric on the contact distribution, a sub-Riemannian manifold. The nilpotentization  of the contact  distribution in a contact manifold of dimension $3$ is   a Heisenberg group $H(1)$. 
  
 The horizontal linear maps on the Heisenberg group are known. Moreover, the group of isometries 
 of $H(1))$ which preserve the origin is   $SU(1) = SO(2)$. 
  
  In order to classify  all homogeneous contact 3 manifolds, we have to consider two cases. 
  The first case is $\mathfrak{g}$ 4 dimensional, with a basis  $\left\{X_{0}, X_{1}, X_{2}, X_{3}\right\}$, such that $X_{0}$ is a basis for $D_{0}$, $\left\{ X_{1}, X_{2}\right\}$ a basis for D, $X_{3}$ a basis 
  for $V_{2}$.  The vectorspace $\mathfrak{g}$ has the direct sum decomposition 
  $$\mathfrak{g} \ = \ D_{0} + D + V_{2}$$
  We shall note $V_{1} = D_{0} + D$. For any $x \in \mathfrak{g}$ we shall use the decomposition 
  $$x = x_{0} + x_{1} + x_{2}$$
  which comes from the direct sum decomposition of $\mathfrak{g}$. 
  
  We also have a Lie bracket $[\cdot, \cdot ]$ on $\mathfrak{g}$ and an Euclidean metric on $D$, called $g$. 
  
  The definition \ref{defhoms} applied in this particular example gives the following relations. 
  There is $a \not = 0$ such that 
   $$[X_{0}, X_{1}] =  a X_{2}$$
   $$[X_{0}, X_{2}] =  -a  X_{1}$$
The dilatations $\delta_{\varepsilon}$ are defined for any $\varepsilon > 0$ and for any $x \in \mathfrak{g}$ by 
$$\delta_{\varepsilon} x \ = \ \varepsilon x_{0} + \varepsilon  x_{1} + \varepsilon^{2} x_{2}$$
Condition (homs - f) implies that 
$$[X_{0}, X_{3}] \ = \ b_{03} X_{0} + e_{03}X_{3}$$
The other bracket relations are: 
$$[X_{1}, X_{2}] \ = \ b_{12} X_{0} + c_{12} X_{1} + d_{12} X_{2} + e_{12} X_{3}$$
$$[X_{2}, X_{3}] \ = \ b_{23} X_{0} + c_{23} X_{1} + d_{23} X_{2} + e_{23} X_{3}$$
$$[X_{1}, X_{3}] \ = \ b_{13} X_{0} + c_{13} X_{1} + d_{13} X_{2} + e_{13} X_{3}$$
Condition (homs - e) implies $e_{12} \not = 0$. 

We have to check the Jacobi relations for the bracket. After a careful computation we are left with the 
following bracket relations: 
$$[X_{0}, X_{1}] =  a X_{2}$$
$$[X_{0}, X_{2}] =  -a  X_{1}$$
$$[X_{0}, X_{3}] = 0$$
$$[X_{1}, X_{2}] = b_{12} X_{0} + e_{12} X_{3}$$
$$[X_{1}, X_{3}] = d X_{2}$$
$$[X_{2}, X_{3}] = -d X_{1}$$
We can modify $X_{0}$ and $X_{3}$ such  such that $a = 1$. This is not modifying the dilatations and the metric $g$. We can then diagonalise the metric $g$. This is leaving us with 2 strictly positive parameters $\lambda_{1}, \lambda_{2}$ associated to the metric 
$$g = \lambda_{1} X_{1}^{*} \otimes X_{1}^{*} + \lambda_{2} X_{2}^{*} \otimes X_{2}^{*}$$
and the parameters $b_{12}$,$d$, $e_{12}$ from the bracket. By rescaling $X_{1}$, $X_{2}$, $X_{3}$ with arbitrary (but not null) numbers $\alpha_{1}, \alpha_{1}, \alpha_{2}$, we see that the set of parameters $$(\lambda_{1}, \lambda_{2}, b_{12},d, e_{12}) \ \ \ \ \ \mbox{ and }  \ \ \ \  (\alpha_{1}^{2} \lambda_{1}, \alpha_{1}^{2} \lambda_{2},  \alpha_{1}^{2} b_{12},  \alpha_{2} d, \frac{\alpha_{1}^{2}}{\alpha_{2}} e_{12})$$
are equivalent. 

The parameter $b_{12}$ will not count when the metric factorisation with respect to the action 
of $SU(1)$ is made, therefore we are left with two parameters: 
$$\frac{d e_{12}}{\lambda_{1}} \ \ \ \ , \ \ \ \ \frac{\lambda_{2}}{\lambda_{1}}$$ 
or equivalently with a point $[d e_{12}, \lambda_{2}, \lambda_{1}]$ from a convex set in projective plane, 
with a free choice for the sign of $d  e_{12}$. 

But this is not all. Recall that in the end we are looking to metric profiles associated to dilatation structure. When any metric profile is defined, a further normalisation takes place: we have to decide what the parametrisation of the metric profile is, namely to choose $\mathbb{P}(1)$. 

This can always be translated, for the metric profile of the metric space $X_{\sigma}$, by a conformal modification of the metric $g$ (see for this relation (\ref{compa3})). If we use the dilatation metric profile 
then what we do is a renormalisation in the sense of relation (\ref{compa2})

This is leaving us with only a one dimensional family of "good" curvatures in this case. 
An explanation of this final normalisations is found in section \ref{whatcurv}. 

The contact homogeneous spaces  we have just discussed correspond to spaces with maximal symmetry, using the terminology of 
Montgomery  \cite{montgome} chapter 7, section 7.10, paragraph "Examples with maximal symmetry".

The second case is $\mathfrak{g}$ 3 dimensional, Lie algebra with the bracket $[\cdot, \cdot]$, and 
a metric $g$. We pick a basis $\left\{ X_{1}, X_{2}, X_{3} \right\}$ and a decomposition 
$$\mathfrak{g} \ = \ D + V_{2}$$
such that $\left\{X_{1}, X_{2}\right\}$  is a basis for $D$ and $X_{3}$ is a basis for $V_{2}$. The dilatations  $\delta_{\varepsilon}$ are defined for any $\varepsilon > 0$ and for any $x \in \mathfrak{g}$ by 
$$\delta_{\varepsilon} x \ = \  \varepsilon  x_{1} + \varepsilon^{2} x_{2}$$
The only constraint on the bracket relations is $$[X_{1}, X_{2}]_{3} \not = 0$$

The classification uses the structure group $G_{0}$ as in Cartan method. We shall not reproduce here the details, but only give the answer. 

Any homogeneous 3 dimensional Lie algebra is isomorphic with a triple $\sigma$ which has the 
following form. 

There is a basis $\left\{ X_{1}, X_{2}, X_{3} \right\}$ such that the metric $g$ has the form 
$$g = X_{1}^{*} \otimes X_{1}^{*}$$ 
the dilatations are given for any $\varepsilon > 0$ by
$$\delta_{\varepsilon} X_{1} \  = \  \varepsilon X_{1} \ \ , \ \ \delta_{\varepsilon} X_{2} \ = \ \varepsilon X_{2} \ \ , \ \  \delta_{\varepsilon} X_{3} \ = \ \varepsilon^{2} X_{3}$$ 
and the bracket relations are
$$[X_{1}, X_{2}] \ = \ X_{3}$$
$$[X_{2}, X_{3}] \ = \ \rho \cos^{2} \phi X_{1} \ + \ \rho \sin \phi \cos \phi X_{2} \ + \ \gamma \cos \phi X_{3}$$
  $$[X_{3}, X_{1}] \ = \ \rho \sin \phi \cos \phi X_{1} \ + \ \rho \sin^{2} \phi  X_{2} \ + \ \gamma \sin \phi X_{3}$$
  where $\rho, \phi, \gamma \in \mathbb{R}$ are arbitrary numbers. 

As in the previous case we have a last normalisation to do, explained in  section \ref{whatcurv}, relations  (\ref{compa2}) and (\ref{compa3}). A straightforward way to  understand this normalisation consists in the use of  the deformed bracket from (homs -e), which here takes the form 
$$[X_{1}, X_{2}]_{\varepsilon} \ = \ X_{3}$$
$$[X_{2}, X_{3}]_{\varepsilon} \ = \ \rho \varepsilon^{2} \cos^{2} \phi X_{1} \ + \ \rho \varepsilon^{2} \sin \phi \cos \phi X_{2} \ + \ \gamma \varepsilon \cos \phi X_{3}$$
  $$[X_{3}, X_{1}]_{\varepsilon} \ = \ \rho \varepsilon^{2} \sin \phi \cos \phi X_{1} \ + \ \rho \varepsilon^{2} \sin^{2} \phi  X_{2} \ + \ \gamma \varepsilon \sin \phi X_{3}$$
Denote by $\sigma_{\varepsilon}$ the homogeneous space with the bracket defined by $[\cdot, \cdot]_{\varepsilon}$, the canonical  metric $g$ and dilatations. This space is described by the parameters $(\rho \varepsilon^{2}, \phi, \gamma \varepsilon)$.  We shall finally identify these spaces, which leaves us with only a 2 dimensional parameter space, again in total agreement with Montgomery, Hughen, Agrachev et. al. {\it op.cit.}.  

Particular examples are:  SO(3), SL(2,R), E(1,1), each with naturally chosen generating 2 dimensional distributions. 

A similar  classification can be done, in an easier way,  using homogeneous ensembles. These are described in next section.

\section{Homogeneous ensembles and deformations}
\label{homes}

In this section we shall describe first the homogeneous ensembles, then  the deformations of such objects.

The inspiration for the notion of a homogeneous ensemble comes from  the construction of a normal frame in a sub-Riemannian manifold. We recall it here. 

We start with a system of vectorfields $\left\{ X_{1}, ... , X_{p} \right\}$ which span the distribution $D$. 
There is also an Euclidean metric $g$ on the distribution $D = V^{1}$. 

We add, in lexicographic order, brackets of the initial vectorfields, $[X_{i}, X_{j}]$, until we obtain a 
basis for $V^{2}$. We can equally extend the metric $g$ to $V^{2}$ by 
$$g([X_{i}, X_{j}], X_{k}) = 0$$
$$g([X_{i}, X_{j}], [X_{i}, X_{j}]) \ = \ g(X_{i}, X_{i}) g(X_{j}, X_{j})$$
if $[X_{i}, X_{j}]$ was added to the normal frame under construction. Moreover, any two  different  vectorfields added at same step of the construction are orthogonal with respect to $g$. 

We repeat the procedure until we complete the normal frame. What we get? 

We end with a frame $\left\{ X_{1}, ... , X_{n} \right\}$ which constitutes the set of nodes of a  
tree with leaves $\left\{ X_{1}, ... , X_{p} \right\}$ and roots $\left\{ X_{i_{1}}, ... , \right\}$ (the vectorfields 
added at the final step of the construction). To any node $X_{k}$ which is not a leaf corresponds two 
branches pointing to $X_{k}^{(1)}$, $X_{k}^{(2)}$,  such that $X_{k}^{(1)} \in \left\{ X_{1}, ... , X_{p} \right\}$ and 
$$X_{k} \ = \ [X_{k}^{(1)}, X_{k}^{(2)}]$$
meaning that $X_{k}$ has been obtained as the bracket of $X_{k}^{(1)}$, $X_{k}^{(2)}$ at some step 
of the construction. 

We also have a degree function associated to the tree in a natural way (distance from the leaves, plus one). 

We can construct  a metric $g$ which extends the metric (denoted with same letter) on the distribution.   The extended metric is unquely defined by the following conditions: 
\begin{enumerate}
\item[-] if $deg \ X_{i} \ \not = \ deg \ X_{j}$ or $deg \ X_{i} \ = \ deg X_{j}$  but $i \not = j$,  then $g(X_{i}, X_{j}) = 0$ 
\item[-] if $deg \ X_{k} \ \geq 2$ then 
$$g(X_{k}, X_{k}) \ = \ g(X_{k}^{(1)}, X_{k}^{(1)}) g(X_{k}^{(2)}, X_{k}^{(2)})$$
\end{enumerate}

The algebraic counterpart of this construction is described further. 
\begin{enumerate}
\item[(home - a)] We have a vectorspace $\mathfrak{g}$ of dimension $n$, endowed with a bilinear form 
$$[ \cdot , \cdot ] : \mathfrak{g} \times \mathfrak{g} \rightarrow \mathfrak{g}$$
which is antisymmetric: for any $u,v \in \mathfrak{g}$ we have 
$$ [v,u] \ = \ - [u,v]$$
\item[(home - b)] The vectorspace $\mathfrak{g}$ admits a gradation 
$$\mathfrak{g} \ = \ V_{0} +  V_{1} + ... + V_{m}$$ 
such that for any $i \in \left\{ 0,1\right\}$ and $j \in \left\{ 0, ... ,m \right\}$,  if $i+j \leq m$ then 
$$V_{0} + ... + V_{i+j} \ = \   V_{0} + ... + V_{j} + [V_{i}, V_{j} ]$$ 
\item[(home - c)] Moreover, for any $i,j, k \in \left\{0, ... ,m \right\}$, if two of the following 3 relations holds 
$$i + j \leq m \ \ , \ \ j+k \leq m \ \ , \ \ k+i \leq m$$
then for any $u \in V_{i}$, $v \in V_{j}$, $w \in V_{k}$  the Jacobi relation is true: 
$$[[u,v],w] \ + \ [[w,u],v] \ + \ [[v,w],u] \ = \ 0$$
\end{enumerate}

We can quickly deduce from (home - b) that   for any $u,v \in V_{0}$ we have $[u,v] \in V_{0}$. Moreover $V_{0}$ is a Lie algebra 
with respect to the bracket $[\cdot , \cdot]$, because the Jacobi identity is satisfied, due to (home -c).

\begin{enumerate}
\item[(home - d)] Denote by $p$ the dimension of $V_{1}$.  The function $u_{0} \in V_{0} \mapsto [u_{0}, \cdot] : V_{1} \rightarrow \mathfrak{g}$ is 
the image of an action of a subalgebra of $\mathfrak{so}(p)$. This means that there is an injective  
algebra morphism  $Q: V_{0} \mapsto \mathfrak{so}(p)$ such that for any $u_{0} \in V_{0}$ and 
$u \in V_{1}$ we have 
$$[u_{0}, u] \ = \ Q(u_{0})u$$
\item[(home - e)] Associated to the gradation of $\mathfrak{g}$ there is a degree map $deg \ : \mathfrak{g} \rightarrow 
\left\{ 0, ... , m \right\}$, defined by: if $u \in V_{0} + V_{1}$ then $deg \ u = 1$; for $k \geq 2$ we have 
$deg \ u \ = \ k$ if and only if  $u \in V_{1} + ... + V_{k}$ but $u \not \in V_{1} + ... + V_{k-1}$. 

This degree map induces a dilatation flow. For any $\varepsilon > 0$ the dilatation 
$\delta_{\varepsilon}$ is a linear map defined by the values it takes on each $V_{k}$, $k = 0, ... , m$: 
if $x \in V_{k}$ then $\delta_{\varepsilon} x \ = \ \varepsilon^{deg \ x} x$. Remark that due to the definition 
of the degree, for any $x \in V_{0}$ we have $deg \ x \ =1$, hence $\delta_{\varepsilon} x = \varepsilon x$. 
\item[(home - f)] We also have a metric $g$ on $\mathfrak{g}$ which is null on $V_{0}$ and stricly positive definite on $V_{1}$. Thi smetric can be extended to a metric which  makes the gradation 
of $\mathfrak{g}$ orthogonal. 
\end{enumerate}

\begin{defi}
A homogeneous ensemble is a triple $([\cdot, \cdot], \delta, g)$ which satisfies all the conditions 
(home - a), ... , (home - f).
\label{defhome}
\end{defi}

We shall explain now what is the pointed metric space associated to a homogeneous ensemble. 

To a triple $\sigma = (\mathfrak{g}, \delta, g)$ we associate a length measure given by the metric $g$,  
and the bracket $[\cdot, \cdot]$. Even if this bracket is not a Lie bracket, we can ignore this and define the distribution at $x \in \mathfrak{g}$ by the formula 
$$D_{x}  \ = \ \left\{ DL^{\sigma}_{x}(0) Y \mbox{ :  } Y \in V_{0} + V_{1} \right\}$$
where the operator $DL^{\sigma}_{x}(0)$ mimicks the derivative of a left translation at $x$. This operator 
has the expression 
$$DL^{\sigma}_{x} (0) Y \ = \ \sum_{k =0}^{\infty} \frac{[x, \cdot]^{k}}{(k+1)!} Y$$ 
which is inspired by the Baker-Campbell-Hausdorff formula. 

This length measure induces a (pseudo) distance denoted by $d$. For any $\varepsilon > 0$ define the (pseudo) distance 
$$d_{\varepsilon}(x,y) \ = \ d (\delta_{\varepsilon}^{-1} x , \delta_{\varepsilon}^{-1} y )$$
Make the metric factorisation with respect to $d$  in a domain of convergence of the operators 
$DL^{\sigma}$. The resulting metric space is (in a neighbourhood of $0$)  a regular sub-Riemannian manifold. 

Declare $a>0$ admissible if the open ball centered at $0$ and radius $1$, defined with respect to the 
distance $d_{2a}$, is contained in (the factor of) the domain of convergence just mentioned.  
Consider now the number 
$$R = \min \left( 1, \sup \left\{a>0  \ \mbox{ admissible } \right\} \right)$$  
Associate then to $\sigma$ the compact metric space $(X_{\sigma}, d_{\sigma})$ where $X_{\sigma}$ is 
the closed ball of radius $1$ centered at $0$ and $d_{\sigma} = d_{R}$.  The (isometry class of ) pointed 
metric space associated to $\sigma$ is $[X_{\sigma}, d_{\sigma}, 0]$. 

The deformation (metric profile) associated to $\sigma$ is the deformation by dilatations based at $0$, according to formula (\ref{1defor}). We shall denote this metric profile by 
$a \mapsto \mathbb{P}(\sigma) (a)$.

\section{What curvature is}
\label{whatcurv}

In preparation for a refined definition of curvature we introduce an action of $GL(\mathfrak{g})$ on 
$\sigma \ = \ ([\cdot , \cdot], \delta, g)$. This is the most natural one, namely the transport. 
Fo any $F \in GL(\mathfrak{g})$ and any  homogeneous ensemble $\sigma$ over $\mathfrak{g}$  we define 
$$F \sigma \ = \ ( F [ F^{-1} \cdot , F^{-1} \cdot] , F \delta F^{-1} , g(F^{-1} \cdot , F^{-1} \cdot ))$$

There is also an action by dilatations, which is defined differently, by 
$$\delta_{\varepsilon} * \sigma \ = \  ( F [ F^{-1} \cdot , F^{-1} \cdot] , F \delta F^{-1} , g)$$
Finally we have an action of $(0,+\infty)$ on homogeneous ensembles, defined by: 
$$\varepsilon . \sigma \ = \ (  [  \cdot , \cdot] ,  \delta , \varepsilon^{2} g)$$
By straightforward computation we have:
\begin{prop}
 For any  $a > 0$ and small enough $\varepsilon > 0$ 
\begin{equation} 
\mathbb{P}(\sigma)(a \varepsilon) \ = \ \mathbb{P}( \delta_{a}^{-1} * \sigma) (\varepsilon)
\label{compa2}
\end{equation}
\begin{equation}
\mathbb{P}^{m}(\sigma)(a \varepsilon) \ = \ \mathbb{P}^{m}(a^{-1} . \sigma) (\varepsilon)
\label{compa3}
\end{equation}
\end{prop}

This shows that the dilatation metric profile $a \mapsto \mathbb{P}(\sigma) (a)$ lies in the image of 
homogeneous ensembles in CMS. It also makes  clear that the final, metric,  normalisation, invoked  in subsections  \ref{rie} and \ref{srcon}, is about the choice of initial ($a=1$) point on the dilatation metric profile of the homogeneous ensemble. 

To the homogeneous ensemble $\sigma$ we can also associate the metric profile in $0 \in \mathfrak{g}$ and we shall denote this by $a \mapsto \mathbb{P}^{m}(\sigma)(a)$. The same remarks as previous hold for the metric profile and the action of $(0,+\infty)$.

We can define several types of curvatures. We shall list here only the {\it metric} curvature and 
the {\it dilatation} curvature. 

\begin{defi}
Let $(X,d)$ be a metric space and $p \in X$ a point such that the metric profile associated to $(X,d)$ 
and $p$ can be prolonged at $\varepsilon = 0$ and it is rectifiable at $\varepsilon = 0$. We shall say 
that the homogeneous ensemble $\sigma$ represents the metric  curvature of $(X,d)$ at $p$  if the metric profile of $(X,d)$ at $p$ is equivalent with the metric profile $\mathbb{P}^{m}(\sigma)$.
\label{defmetric} 
\end{defi}
 
 \begin{defi}
Let $(X,d)$ be a metric space endowed with a dilatation structure and $p \in X$ a point such that the dilatation metric profile associated to $(X,d)$, the dilatation structure  
and $p$ can be prolonged at $\varepsilon = 0$ and it is rectifiable at $\varepsilon = 0$. We shall say 
that the homogeneous ensemble $\sigma$ represents the dilatation curvature of $(X,d)$ at $p$  if the dilatation  metric profile of $(X,d)$ at $p$ is equivalent with the dilatation metric profile $\mathbb{P}(\sigma)$.
\label{defdila} 
\end{defi}

In the case of Riemannian manifolds, a basic lemma asserts the existence of special coordinates systems which are adapted around an open open neighbourhood of an arbitrary point $p$  from  the manifold. Moreover, this coordinated can be chosen such that  at $p$ are adapted "up to order 2" 
(not defined in this paper but straightforward). This guarantees that the metric and dilatation curvature 
coincide in the sense that the same Riemannian homogeneous ensemble represents the metric and the 
dilatation curvature. We doubt very much that this is possible for general sub-Riemannian manifolds. 

Remark also that in these definitions of curvature there is contained information about the tangent space. In the Riemannian case this is superfluous but for sub-Riemannian manifolds it is not.

 For example the Heisenberg group, which is Carnot group, is itself a representant of its curvature. In a sense it has curvature $0$ (the dilatation and the metric profile are both stationary), but it is different from an Euclidean plane which has also $0$ curvature.

In the case of a Riemannian manifold with singularities, the metric profile at a singularity point is an 
Euclidean cone. Form the point of view of this paper Euclidean cones have stationary metric profiles 
(in some sense curvature $0$), but different from the metric profile of an Euclidean plane. 

 It is useful to picture the curvature as a vector tangent to the metric or dilatation profile: the orientation and length of the curvature "vector" corresponds to second 
order infinitesimal informations about the pointed metric space. The base point of the "vector"  contains 
first order informations, namely the description of the tangent space. 

The problem of computation of dilatation curvature of a regular sub-Riemannian manifold is left for 
future work. Some things are however clear: such computation should involve a normal frame and the 
dilatation structure associated to it. The real problem involves the fact that the bracket associated with 
a normal frame has nonconstant coefficients. 

It is however important to notice that eventual differential geometric   notions of curvature connected to  normal frames  and privileged coordinates are in fact expressing the dilatation --- not the metric --- curvature. 
 
\section{Curvature and coadjoint orbits} 
 \label{lastsec}

 In this section we show that classification of dilatation curvatures means classification of some coadjoint orbits representations. This points directlly to the powerful orbit method and constitute a second link between sub-Riemannian geometry and quantum mechanics. The first link can be uncovered from 
 Buliga \cite{buliga2} section 5 "Case of the Heisenberg group". The results from the mentioned section can be easily generalized for a pre-quantum contact manifold and are the subject of a paper in preparation.

 To say it in few words: a quantum dynamical system is just a dynamical system  $t \mapsto \phi_{t} : X \rightarrow X$ in a metric measure space $(X,d, \mu)$ endowed with a dilatation structure $\delta$. The dynamical system has to be:
 \begin{enumerate}
 \item[(A)] measure preserving, 
 \item[(B)] smooth (orbits and the transformations $\phi_{t}$ should be $\delta$ derivable), 
 \item[(C)]  the orbits 
 $$\left\{ \phi_{t}(x) \mbox{ :  } t \in [a,b]\right\}$$
  have Hausdorff dimension 2. The Hausdorff measure 2 should be absolutely continuous with 
 respect to $t$ and the density of this measure is by definition the Hamiltonian. 
 \end{enumerate}
 
 Look for example to  the case of a $S^{1}$-bundle associated to a pre-quantum contact manifold. 
 We see see this as a contact sub-Riemannian manifold. 
 
  If the pre-quantum contact manifold is a pre-quantization of an integral symplectic manifold then  any  dynamical system which satisfies (A), (B),   is a lifting  of a Hamiltonian dynamical system on the 
 symplectic manifold. Moreover, in this case the condition (C) is satisfied, in the sense that the orbits of the dynamical system have indeed Hausdorff dimension 2 and the density of the Hausdorff measure 
 $2$ with respect to the (transport from $\mathbb{R}$ of the) Lebesgue measure $1$ on the curve is a Hamiltonian for the  dynamical system on the symplectic manifold.
 
 In the case of a Riemannian manifold such dynamical systems  correspond to random walks or 
 random dynamical systems. 
 
A measurement process should correspond to trying to make an Euclidean chart of this dynamical 
system. If the space is not Euclidean at any scale (as the sub-Riemannian Heisenberg group, for example) then such a map is impossible to be done exactly (i.e. in a derivable way). 

Planck constant might be the effect of this fact, namely it could measure the distance from the (metric profile or dilatation profile) and best Euclidean approximations. After reading this paper one can be sensible to the idea that the Planck constant could measure a distance between curvatures.

Let us come back to the subject of this paper. We want to associate to a homogeneous ensemble a 
coadjoint orbit representation. Let us imagine that the homogeneous ensemble corresponds to a normal frame associated to a regular sub-Riemannian manifold. We can change the normal frame without changing the dilatation structure. Indeed, it is enough to change the basis of the distribution and then 
build another normal frame with the help of the new basis of the distribution. 

This corresponds in terms of homogeneous ensembles to the action on the homogeneous ensemble 
$\sigma = ([\cdot, \cdot], \delta , g)$ of a group $G(\sigma) \subset GL(\mathfrak{g})$ such that the profile 
$\mathbb{P}(\sigma)$ is left unchanged.

Consider the decomposition of $\mathfrak{g}$: 
$$\mathfrak{g} \ = \ V_{0} + V_{1} + .... + V_{m}$$
and the algebra morphism $Q : V_{0} \rightarrow so(p)$ such that for any 
$u_{0} \in V_{0}$ and $u_{1} \in V_{1}$ we have 
$$[u_{0}, u_{1}] \ = \ Q(u_{0}) u_{1}$$
Finally,  consider a basis $\left\{ X_{1}, ... , X_{n} \right\}$ of $V_{1} + ... + V_{m}$, adapted to the gradation (i.e. there is a partition of the basis in $m$ sets, each forming a basis for one of the $V_{i}$). This basis is chosen such that 
it  satisfies all the conditions explained previously in the construction of a normal frame, at the beginning 
of section  \ref{homes}. 

The basis  $\left\{ X_{1}, ... , X_{n} \right\}$  constitutes the set of nodes of a  
tree with leaves the basis $\left\{ X_{1}, ... , X_{p} \right\}$ of $V_{1}$  and roots $\left\{ X_{i_{1}}, ... , \right\}$ (the basis of $V_{m}$). To any node $X_{k}$ which is not a leaf corresponds two 
branches pointing to $X_{k}^{(1)}$, $X_{k}^{(2)}$,  such that $X_{k}^{(1)} \in \left\{ X_{1}, ... , X_{p} \right\}$ and 
$$X_{k} \ = \ [X_{k}^{(1)}, X_{k}^{(2)}]$$
meaning that $X_{k}$ has been obtained as the bracket of $X_{k}^{(1)}$, $X_{k}^{(2)}$ at some step 
of the construction. 

We also have a degree map associated to the tree in a natural way (distance from the leaves, plus one). 
The degree map of the tree coincide with 
the degree map of the gradation $V_{1} + ... + V_{m}$.  Extend the metric $g$ on $V_{1} + ... V_{m}$ by  the metric relation 
$$g(X_{k}, X_{k}) \ = \ g(X_{k}^{(1)}, X_{k}^{(1)}) g(X_{k}^{(2)}, X_{k}^{(2)})$$
and by the condition that the gradation is a $g$ orthogonal decomposition. 

 The metric $g$ has the expression 
$$g \ = \ \sum_{k=1}^{m} \sum_{deg \ X_{i}  = \ deg \ X_{j}  = k} g^{k}_{ij} X_{i}^{*} \otimes X_{j}^{*}$$
We shall add to this metric the minus Killing form on $\mathfrak{so(p)}$, transported back by the 
morphism $Q$ on $V_{0}$. We get a strictly positive definite metric $\bar{g}$. 

We are going to define now the group $G(\sigma)$. An invertible linear transformation $F : \mathfrak{g} \rightarrow \mathfrak{g}$ belongs to $G$ if and only if: 
\begin{enumerate}
\item[(a)] for any $k = 0, ... , m$ we have 
$$F(V_{0} + ... + V_{k}) \ = \  V_{0} + ... + V_{k}$$
\item[(b)] for any $u_{0} \in V_{0}$ we have $F(u_{0}) = u_{0}$.
\item[(c)] for any $u_{0} \in V_{0}$ and any $u_{1} \in V_{1}$ we have 
$$F[u_{0}, u_{1}] \ = \ [u_{0}, F(u_{1})]$$
that is $F$ commutes with the representation $Q$. This implies that $F(V_{1}) \ = \ V_{1}$. 
\item[(d)] the restriction of $F$ on $ V_{1}$ is a $g$ isometry. 
\end{enumerate}

In contrast with relations (\ref{compa2}) and (\ref{compa3}) we have  
\begin{prop}
For any $F \in G(\sigma)$ and any sufficiently small $\varepsilon > 0$ we have $\mathbb{P}(F\sigma) (\varepsilon) \ = \ \mathbb{P}(\sigma)(\varepsilon)$. 
\end{prop}
The proof is a straightforward computation. 

Let us remark that $G(\sigma) \ = G(\sigma_{\varepsilon})$ for any $\varepsilon > 0$.

To a dilatation curvature  $\mathbb{P}(\sigma)$ we associate in a bijective way the function  
$$ \varepsilon > 0 \ \mapsto \ \mathcal{O}(\sigma_{\varepsilon}) \ = \ \left\{ F \sigma_{\varepsilon} \mbox{ :  } F \in G(\sigma) \right\}$$
where 
$$\sigma_{\varepsilon} \ = \ \delta_{\varepsilon}^{-1}* \sigma$$
We shall call this function a scaled orbit and we denote it by $\mathcal{O}(\sigma)$. 

The purpose of this section is to prove the following theorem. 

\begin{thm}
The action of $G(\sigma)$ on the scaled orbit $\mathcal{O}(\sigma)$ is a coadjoint orbit action. 
\label{tcoad}
\end{thm}

In the following we shall prove the theorem. Fix the  Euclidean metric $\bar{g}$  on $\mathfrak{g}$ and denote it  further by $\left( \cdot , \cdot \right)$. We shall use a basis $\left\{ X_{1}, ... , X_{n}\right\}$, as 
explained previously, to do computations.

Let us  use the notation $\delta_{\varepsilon}^{-1} * \sigma \ =  \ ([ \cdot, \cdot]_{\varepsilon}, \delta, g)$. The dilatations commute, therefore they not change when $\delta_{\varepsilon}$ is applied.

To the bracket $[\cdot , \cdot ]_{\varepsilon}$ associate the linear map $W_{\epsilon}: \mathfrak{g} \rightarrow \mathfrak{gl}(\mathfrak{g})$ by the usual procedure: for any $u, x,y \in \mathfrak{g}$ we have 
$$\left(u, [x,y]_{\varepsilon} \right) \ = \ \left( W_{\varepsilon}(x) u, y\right)$$
In case $\sigma$ is a homogeneous space the map $-W_{\varepsilon}$ is a Lie algebra morphism, as a consequence of the Jacobi identity. In the general case the function $-W_{\varepsilon}$ does have only limited morphism properties, due to (home - c).

 It is comfortable to consider $\varepsilon$ as a variable and $W_{\varepsilon}(x)$ as a polynomial in $\varepsilon$. This is indeed a polynomial,  because  due to conditions (home - b) and  (home -f) we can prove that the limit exists  
 \begin{equation}
 \lim_{\varepsilon \rightarrow 0} [u,v]_{\varepsilon}  \ = \ [u,v]_{N}
 \label{signilp}
 \end{equation}
 This is same relation as the one in (homs - e). See also   relation (\ref{compa2}). We state this as a proposition. 

\begin{prop}
The function which maps $u \in \mathfrak{g}$ to  $ W(x) \in gl(\mathfrak{g})[\varepsilon]$ is well defined. 
\end{prop}
The function $W$ encodes all the information about the dilatation curvature induced by $\sigma$. 

This function   can be seen as a linear space 
$$\mathcal{B}(\sigma) \   \subset \  \left( gl(\mathfrak{g}[\varepsilon]) \times \mathfrak{g} \right)^{*}$$  
in the dual of the natural semidirect product of $gl(\mathfrak{g}[\varepsilon])$ with $\mathfrak{g}$, defined by  
$$ \mathcal{B} (\sigma) \ = \ \left\{ \left(  \begin{array}{cc} W_{\varepsilon}(u) & 0 \\ 
u & 0 \end{array} \right) \mbox{ : } u  \in \mathfrak{g} \right\}$$ 

We call $\mathcal{B}(\sigma)$ the {\it bunch} associated to $\sigma$. 
 
It is easy to see that the action of $G(\sigma)$ on $\mathcal{O}(\sigma)$ transforms in the coadjoint action of 
$$G(\sigma) \ \equiv \ \left\{ \tilde{F} \ = \ \left( \begin{array}{cc} F^{T} & 0 \\ 
0 & 1 \end{array} \right) \mbox{ : } F \in G(\sigma) \right\}$$
namely we have the relation 
$$Ad^{*}_{\tilde{F}} \mathcal{B}(\mu) \ = \ \mathcal{B}(F \mu)$$
for any $\mu \in \mathcal{O}(\sigma)$ and any $F \in G(\sigma)$. The proof of the theorem is finished.

The {\it whole bunch} of $\sigma$ is by definition 
$$\mathcal{WB}(\sigma) \ = \ \bigcup_{\mu \in \mathcal{O}(\sigma)} \mathcal{B}(\sigma)$$
It is a collection of (scaled) coadjoint orbits and therefore it has good chances to be 
 a Poisson $G(\sigma)$-manifold.

To any coadjoint orbit corresponds a natural  representation. The orbit method of Kirillov is 
a guide towards classification of unitary representations using coadjoint orbits. See for this 
Kirillov \cite{kir1} section 15, using as a guide Kirillov \cite{kir2} sections 1 and 2. For the induced representation notion see Kirillov \cite{kir1}, section 13. There is a huge literature dedicated to the orbit 
method. We want here just to make the connection between curvature of sub-Riemannian spaces and 
some coadjoint orbit representations, with the hope what the link will pay back in the future, by translating techniques and objects related to the orbit method to the domain of sub-Riemanian geometry.  

Let us see how a  scaled orbit induces a representation by the prequantization technique. 

The moment map associated to the action of $G(\sigma)$ on the whole bunch $\mathcal{WB}(\sigma)$ is given by the inclusion of the whole bunch in $\left( gl(\mathfrak{g}[\varepsilon]) \times \mathfrak{g} \right)^{*}$. We identify this dual with the Lie algebra by he natural metric that we have, namely the 
usual one induced by trace  on $\mathfrak{gl}(\mathfrak{g})$ and the metric $\left( \cdot , \cdot \right)$ on 
$\mathfrak{g}$. We shall denote $(A,B) \ = \ \tr(AB^{T})$ for any $A,B \in \mathfrak{gl}(\mathfrak{g})$. 

Let $(W_{\varepsilon}, u) \in \mathcal{WB}(\sigma)$. The moment map 
$$J : \mathcal{WB}(\sigma) \rightarrow \ \left(Lie \ G(\sigma) \right)^{*}[\varepsilon]$$
has the following value at $(W_{\varepsilon}, u)$: 
$$\forall f \in \ Lie \  G(\sigma) \ \ \ \ \langle J \left(W_{\varepsilon}, u\right) , f \rangle\ = \   \left( W_{\varepsilon} ,  f \right)$$
The prequantization associates to any $f \in \ Lie \ G(\sigma)$ a self-adjoint operator $Q(f)$ on $C^{\infty}(\mathcal{WB}(\sigma))$ in the following way: for any $h \in C^{\infty}(\mathcal{WB}(\sigma))$ 
we have 
$$Q(f) h \ = \ \frac{i}{2\pi} \frac{d}{dt \mid}_{t=0} \left( h \circ \exp (t \, f) \right) \ + \ \langle J , f \rangle h$$
Let us  take as the representation space the following: 
$$\mathcal{S}(\sigma) \ = \ \left\{ h_{\mid_{\mathcal{WB}(\sigma)}} \mbox{ : }  h \in C^{\infty} \left(\mathfrak{gl}(\mathfrak{g})[\varepsilon] \times \mathfrak{g} \right)\right\}$$
or a closure in a Hilbert space norm. 

Tthen the operator $Q(f)$ takes the form: 
\begin{equation}
Q(f) h \ =  \  \frac{i}{2 \pi} \left\{  (\frac{\partial h}{\partial W} ,  [ W_{\varepsilon}, f] ) \ + \ 
(\frac{\partial h}{\partial u}, f \ u ) \right\} \ + \  \left( W_{\varepsilon} , f \right) \ h 
\label{quant}
\end{equation}
Actually we may take instead of $\mathfrak{gl}(\mathfrak{g})[\varepsilon]$ only a finite dimensional subspace of polynomials up to a certain degree (for example $2m -1$). 

This is the representation which is associated to the dilatation profile induced by the homogeneous 
ensemble $\sigma$. 

If the homogeneous ensemble is a cone (the dilatations commute with the bracket) then in the representation we do not have dependence on $\varepsilon$.

Let us go back to the group $G(\sigma)$. In fact, by relation (\ref{signilp}) we have  
$$G(\sigma) \ = \ G(\sigma_{N}) \ \ , \ \ \ \sigma_{N} \ =  \  ([\cdot ,\cdot]_{N} , \delta , g)$$
The cone $\sigma_{N}$ is associated with the metric tangent space to the space $X_{\sigma}$, at point 
$0$. 

Let us therefore consider a regular sub-Riemannian manifold with the same (up to isometry!) metric tangent space in any of its points, which corresponds to $\sigma_{N}$. Suppose that, when endowed with a dilatation structure, the manifold admits dilatation curvature in any point. Then the dilatation curvatures are classified by self-adjoint  representations of  $Lie \ G(\sigma_{N})$.  This can be stated 
as a theorem: 

\begin{thm}
The possible dilatation curvatures of a regular sub-Riemannian manifold $M$ at a point $x$ are classified by self-adjoint representations of the group $G(\sigma_{N}(x))$, where $\sigma_{N}(x)$ is the Carnot group which represents the metric tangent space of $M$ at $x$. 
\end{thm}

In the following we mean by "unitary dual" of a group the class of self-adjoint representations of the Lie algebra of the group, factorized by equivalence of representations. We have therefore the following characterization  of dilatation curvatures. 

Consider the class $K(\sigma_{N})$  of all homogeneous ensembles $\sigma$ with a given nilpotentization $\sigma_{N}$, with the action of dilatations by 
$$(0, +\infty) \times K(\sigma_{N}) \rightarrow K(\sigma_{N}) \ \ , \ \ (\varepsilon, \sigma) \ \mapsto \sigma_{\varepsilon}$$

\begin{thm}
The function which associates to any $\sigma \in K(\sigma_{N})$ the self-adjoint representation 
given by (\ref{quant}) with $\varepsilon =1$,  transforms integral curves of the action by dilatations 
into polynomial curves in the unitary dual of $G(\sigma_{N})$. 
\end{thm}

\section{Appendix: Uniform and conical groups}

We start with the following setting: $G$ is a topological group endowed with an uniformity such that the operation is uniformly continuous.  More specifically,
we introduce first the double of $G$, as the group $G^{(2)} \ = \ G \times G$ with operation
$$(x,u) (y,v) \ = \ (xy, y^{-1}uyv)$$
The operation on the group $G$, seen as the function
$$op: G^{(2)} \rightarrow G \ , \ \ op(x,y) \ = \ xy$$
is a group morphism. Also the inclusions:
$$i': G \rightarrow G^{(2)} \ , \ \ i'(x) \ = \ (x,e) $$
$$i": G \rightarrow G^{(2)} \ , \ \ i"(x) \ = \ (x,x^{-1}) $$
are group morphisms.

\begin{defi}
\begin{enumerate}
\item[1.]
$G$ is an uniform group if we have two uniformity structures, on $G$ and
$G^{2}$,  such that $op$, $i'$, $i"$ are uniformly continuous.

\item[2.] A local action of a uniform group $G$ on a uniform  pointed space $(X, x_{0})$ is a function
$\phi \in W \in \mathcal{V}(e)  \mapsto \hat{\phi}: U_{\phi} \in \mathcal{V}(x_{0}) \rightarrow
V_{\phi}  \in \mathcal{V}(x_{0})$ such that:
\begin{enumerate}
\item[(a)] the map $(\phi, x) \mapsto \hat{\phi}(x)$ is uniformly continuous from $G \times X$ (with product uniformity)
to  $X$,
\item[(b)] for any $\phi, \psi \in G$ there is $D \in \mathcal{V}(x_{0})$
such that for any $x \in D$ $\hat{\phi \psi^{-1}}(x)$ and $\hat{\phi}(\hat{\psi}^{-1}(x))$ make sense and   $\hat{\phi \psi^{-1}}(x) = \hat{\phi}(\hat{\psi}^{-1}(x))$.
\end{enumerate}

\item[3.] Finally, a local group is an uniform space $G$ with an operation defined in a neighbourhood of $(e,e) \subset G \times G$ which satisfies the uniform group axioms locally.
\end{enumerate}
\label{dunifg}
\end{defi}
Remark that a local group acts locally at left (and also by conjugation) on itself.

This definition deserves an explanation. 

An uniform group, according to the definition \eqref{dunifg}, is a group $G$ such that left translations are uniformly continuous functions and the left action of $G$ on itself is uniformly continuous too. 
In order to precisely formulate this we need two uniformities: one on $G$ and another on $G \times G$. 

These uniformities should be compatible, which is achieved by saying that $i'$, $i"$ are uniformly continuous. The uniformity of the group operation is achieved by saying that the $op$ morphism is uniformly continuous. 

The particular choice of the operation on $G \times G$ is not essential at this point, but it is 
justified by the case of a Lie group endowed with the CC distance induced by a left invariant distribution. We shall construct a natural CC distance on $G \times G$, which is left invariant 
with respect to the chosen operation on $G \times G$. These distances induce uniformities which transform $G$ into an uniform group according to definition \eqref{dunifg}. 

In  proposition \eqref{opsm} we shall prove  that the operation function $op$ is derivable, even if 
right translations are not "smooth", i.e. commutative smooth according to definition 
\eqref{fdcd}. This will motivate the choice of the operation on $G \times G$. It also gives a hint about what  
a sub-Riemannian Lie group should be.

We prepare now the path to this result. The "infinitesimal version" of an uniform group is a conical 
local uniform group.

\begin{defi}
A conical local uniform group $N$ is a local group with a local action of
$(0,+\infty)$ by morphisms $\delta_{\varepsilon}$ such that
$\displaystyle \lim_{\varepsilon \rightarrow 0} \delta_{\varepsilon} x \ = \ e$ for any
$x$ in a neighbourhood of the neutral element $e$.
\end{defi}

We shall make the following  hypotheses on the local uniform group $G$: there is a local action of $(0, +\infty)$ (denoted by
$\delta$), on $(G, e)$ such that
\begin{enumerate}
\item[H0.] the limit  $\lim_{\varepsilon \rightarrow 0} \delta_{\varepsilon} x \ = \ e$ exists and is uniform with respect to $x$.
\item[H1.] the limit
$$\beta(x,y) \ = \ \lim_{\varepsilon \rightarrow 0} \delta_{\varepsilon}^{-1}
\left((\delta_{\varepsilon}x) (\delta_{\varepsilon}y ) \right)$$
is well defined in a neighbourhood of $e$ and the limit is uniform.
\item[H2.] the following relation holds
$$ \lim_{\varepsilon \rightarrow 0} \delta_{\varepsilon}^{-1}
\left( ( \delta_{\varepsilon}x)^{-1}\right) \ = \ x^{-1}$$
where the limit from the left hand side exists in a neighbourhood of $e$ and is uniform with respect to $x$.
\end{enumerate}

These axioms are the prototype of a dilatation structure. Further comes a proposition which corresponds to theorem \ref{tcon}. 
\begin{prop}
Under the hypotheses H0, H1, H2 $(G,\beta)$ is a conical local uniform group.
\label{here3.4}
\end{prop}

\paragraph{Proof.}
All the uniformity assumptions permit to change at will the order of taking
limits. We shall not insist on this further and we shall concentrate on the
algebraic aspects.

We have to prove the associativity, existence of neutral element, existence of inverse and the property of being conical. The proof is straightforward.
For the associativity $\beta(x,\beta(y,z)) \ = \ \beta(\beta(x,y),z)$ we compute:
$$\beta(x,\beta(y,z)) \ = \ \lim_{\varepsilon \rightarrow 0 , \eta \rightarrow 0} \delta_{\varepsilon}^{-1} \left\{ (\delta_{\varepsilon}x) \delta_{\varepsilon/\eta}\left( (\delta_{\eta}y) (\delta_{\eta} z) \right) \right\}$$
We take $\varepsilon = \eta$ and we get
$$ \beta(x,\beta(y,z)) \ = \ \lim_{\varepsilon \rightarrow 0}\left\{
(\delta_{\varepsilon}x) (\delta_{\varepsilon} y) (\delta_{\varepsilon} z) \right\}$$
In the same way:
$$\beta(\beta(x,y),z) \ = \ \lim_{\varepsilon \rightarrow 0 , \eta \rightarrow 0} \delta_{\varepsilon}^{-1} \left\{ (\delta_{\varepsilon/\eta}x)\left( (\delta_{\eta}x) (\delta_{\eta} y) \right) (\delta_{\varepsilon} z) \right\}$$
and again taking $\varepsilon = \eta$ we obtain
$$\beta(\beta(x,y),z) \ = \  \lim_{\varepsilon \rightarrow 0}\left\{
(\delta_{\varepsilon}x) (\delta_{\varepsilon} y) (\delta_{\varepsilon} z) \right\}$$
The neutral element is $e$, from H0 (first part): $\beta(x,e) \ = \beta(e,x) \ = \ x$. The inverse of $x$ is $x^{-1}$, by a similar argument:
$$\beta(x, x^{-1})  \ = \ \lim_{\varepsilon \rightarrow 0 , \eta \rightarrow 0} \delta_{\varepsilon}^{-1} \left\{ (\delta_{\varepsilon}x)
\left( \delta_{\varepsilon/\eta}(\delta_{\eta}x)^{-1}\right) \right\}$$
and taking $\varepsilon = \eta$ we obtain
$$\beta(x, x^{-1})  \ = \ \lim_{\varepsilon \rightarrow 0}
\delta_{\varepsilon}^{-1} \left( (\delta_{\varepsilon}x) (\delta_{\varepsilon}x)^{-1}\right) \ = \ \lim_{\varepsilon \rightarrow 0} \delta_{\varepsilon}^{-1}(e) \ = \ e$$
Finally, $\beta$ has the property:
$$\beta(\delta_{\eta} x, \delta_{\eta}y) \ = \ \delta_{\eta} \beta(x,y)$$
which comes from the definition of $\beta$ and commutativity of multiplication
in $(0,+\infty)$. This proves that $(G,\beta)$ is conical.
\quad $\blacksquare$

We arrive at a natural realization of the tangent space to the neutral element.
Let us denote by $[f,g] \ = \ f \circ g \circ f^{-1} \circ g^{-1}$ the commutator of two transformations. For the group we shall denote by
$L_{x}^{G} y \ = \ xy$ the left translation and by $L^{N}_{x}y \ = \ \beta(x,y)$. The preceding proposition tells us that $(G,\beta)$ acts locally by left
translations on $G$. We shall call the left translations with respect to the group operation $\beta$ "infinitesimal". Those infinitesimal translations admit
the very important representation:
\begin{equation}
\lim_{\lambda \rightarrow 0} [L_{(\delta_{\lambda}x)^{-1}}^{G}, \delta_{\lambda}^{-1}] \ = \ L^{N}_{x}
\label{firstdef}
\end{equation}

\begin{defi}
The group $VT_{e}G$ formed by all transformations $L_{x}^{N}$ is called the virtual tangent space at $e$ to $G$.
\end{defi}

The virtual tangent space $VT_{x}G$ at $x \in G$ to $G$ is obtained by translating the group operation and the dilatations from $e$ to $x$. This means: define a new operation on $G$ by
$$y \stackrel{x}{\cdot} z \ = \ y x^{-1}z$$
The group $G$ with this operation is isomorphic to $G$ with old operation and
the left translation $L^{G}_{x}y \ = \ xy$ is the isomorphism. The neutral element is $x$.
Introduce also the dilatations based at $x$ by
$$\delta_{\varepsilon}^{x} y \ = \ x \delta_{\varepsilon}(x^{-1}y)$$
Then $G^{x} \ = \ (G,\stackrel{x}{\cdot})$ with the group of dilatations $\delta_{\varepsilon}^{x}$ satisfy the axioms Ho, H1, H2. Define then the virtual tangent
space $VT_{x}G$ to be: $VT_{x}G \ = \ VT_{x} G^{x}$. A short computation shows that
$$VT_{x} G \ = \ \left\{ L^{N,x}_{y} \ = \ L_{x} L^{N}_{x^{-1}y} L_{x} \mbox{ : } y \in U_{x} \in \mathcal{V}(X) \right\}$$
where
$$L^{N,x}_{y} \ = \ \lim_{ \lambda  \rightarrow 0} \delta_{\lambda}^{-1,x} [\delta_{\lambda}^{x}, L_{(\delta_{\lambda}x)^{x, -1}}^{G}] \delta_{\lambda}^{x}$$

We shall introduce the notion of commutative smoothness, which contains
a derivative resembling with Pansu derivative. This definition is a little bit stronger than the one given by Vodopyanov \& Greshnov \cite{vodopis2}, because  their definition is good for a general CC space, when uniformities are taken according to the distances in CC spaces $G^{(2)}$ and $G$.

\begin{defi}
A function $f: G_{1} \rightarrow G_{2}$ is commutative smooth at $x \in G_{1}$, where
$G_{1}, G_{2}$ are two groups satisfying H0, H1, H2,  if the application
$$u \in G_{1} \ \mapsto \ (f(x), Df(x)u) \in G_{2}^{(2)}$$
exists, where
$$Df(x)u \ = \ \lim_{\varepsilon \rightarrow 0} \delta_{\varepsilon}^{-1}
\left(f(x)^{-1}f(x \delta_{\varepsilon}u)\right)$$
and the convergence is uniform with respect to $u$ in compact sets.
\label{fdcd}
\end{defi}

For example the left translations $L_{x}$ are commutative smooth and
the derivative equals identity. If we want to see how the derivative moves
the virtual tangent spaces we have to give a definition.

Inspired by \eqref{firstdef},  we shall introduce the virtual tangent. We proceed as follows: to $f: G \rightarrow G$ and $x \in G$ let associate the function:
$$\hat{f}^{x}: G \times G \rightarrow G \ , \ \ \hat{f}^{x}(y,z) \ = \ \hat{f}^{x}_{y}(z) \ = \ \left(f(x)\right)^{-1}f(xy)z$$
To this function is associated a flow of left translations
$$ \lambda > 0 \ \mapsto \ \hat{f}^{x}_{\delta_{\lambda} y}: G \rightarrow G$$

\begin{defi}
The function $f: G \rightarrow G$ is virtually  derivable at $x \in G$ if there is
a virtual tangent $VDf(x)$ such that
\begin{equation}
\lim_{\lambda \rightarrow 0} \left[ \left(\hat{f}^{x}_{\delta_{\lambda} y} \right)^{-1} , \delta_{\lambda}^{-1} \right] \ = \ VDf(x)y
\label{2nddef}
\end{equation}
\label{vcd}
and the limit is uniform with respect to $y$ in a compact set.
\end{defi}

Remark that in principle the right translations are not commutative smooth. In Buliga \cite{buliga2}, section 4, it is shown that right translations are smooth in the "mild" sense. 

Now that we have a model for the tangent space to $e$ at $G$, we can show that
the operation is commutative smooth.

\begin{prop}
Let $G$ satisfy H0, H1, H2 and $\delta_{\varepsilon}^{(2)} : G^{(2)} \rightarrow G^{(2)}$ be defined by
$$\delta_{\varepsilon}^{(2)} (x,u) \ = \ (\delta_{\varepsilon}x,
\delta_{\varepsilon} y)$$
Then $G^{(2)}$ satisfies H0, H1, H2, the operation ($op$ function) is commutative smooth  and we have the relation:
$$D \ op \ (x,u) (y,v) \ = \ \beta(y,v)$$
\label{opsm}
\end{prop}

\paragraph{Proof.}
It is sufficient to use the morphism property of the operation. Indeed, the right hand side of the relation to be proven is
$$RHS \ = \ \lim_{\varepsilon \rightarrow 0}
\delta_{\varepsilon}^{-1} \left( op(x,u)^{-1} op(x,u) op \left(\delta_{\varepsilon}^{(2)}(y,v)\right)\right) \ = $$
$$=  \ \lim_{\varepsilon \rightarrow 0}
\delta_{\varepsilon}^{-1} \left( op(\delta_{\varepsilon}^{(2)}(y,v))\right) \ = \ \beta(y,v)$$
The rest is trivial.
\quad $\blacksquare$

This proposition justifies the choice of the operation on $G^{(2)} = G \times G$ and it is a quite surprising result. 

We finish this appendix with a question for the reader who consider the material too elementar: consider instead of a (compact) Lie group an uniform group with a dilatation structure and instead of the Lie algebra of the group consider a homogeneous space or ensemble.  What 
modifications to the notion of Hopf algebra should be made in order to recover the duality between 
the universal enveloping algebra of the Lie algebra of a compact Lie group and the commutative Hopf algebra of the group?

\end{document}